\documentclass{article}
\usepackage{amsmath,amsthm,amsfonts,amscd,amssymb,eucal,latexsym}

\begin{document}

\title{Embedding ergodic actions of compact quantum groups on
$C^*$--algebras into  quotient spaces}
\author{Claudia Pinzari\\
Dipartimento di Matematica, Universit\`a di Roma La Sapienza\\
00185--Roma, Italy}
\date{}
\maketitle

\begin{abstract}
The notion of compact quantum subgroup is revisited and an alternative
definition is given. Induced 
representations are considered and a Frobenius reciprocity theorem is
obtained. A relationship between ergodic actions of compact quantum groups 
on $C^*$--algebras and topological transitivity is investigated. A 
sufficient condition for embedding such  actions  in quantum
quotient spaces is obtained.
\end{abstract}

\begin{section}{Introduction}
Consider an ergodic action of  a compact group $G$ by automorphisms
on the  $C^*$--algebra ${\cal C}(X)$
of continuos functions on the compact space $X$.
It is known that this action arises from a transitive right
action of $G$ on $X$ by homeomorphisms. Therefore the stabilizer 
of a point of $X$ is a closed subgroup $K$ of $G$, and 
$X$ can be identified with the right coset space $K\backslash G$ as a
$G$--space.

 The aim of this note is  to 
understand an appropriate generalization of the above property
to ergodic actions of compact quantum groups on
noncommutative $C^*$--algebras. Our interest in this problem arises
from the program of formulating an abstract duality theory for compact
quantum groups 
 where such a 
generalization is needed,
\cite{DPR}, \cite{Pinzari}.

The relationship between topological transitivity and ergodicity in the
case where
$G$ is still a classical compact group but acting on a noncommutative 
$C^*$--algebra ${\cal C}$, has been investigated by Longo and Peligrad in
\cite{LP}.
They proved that ergodicity is  equivalent to the lack of proper closed
$G$--invariant left ideals in ${\cal C}$.

A general theory of ergodic actions of compact matrix pseudogroups
on $C^*$--algebras has been studied by Boca \cite{Boca}, who proved, among
other 
things, that the spectral subspaces of the action are finite dimensional.

The notion of quantum subgroup and quantum quotient space was first
introduced by Podles in
\cite{Podles} for compact matrix pseudogroups, who computed all the
subgroups
and quotient spaces of the quantum $SU(2)$ and $SO(3)$ groups.

Later Wang \cite{Wang}, while  studying ergodic actions of the universal
quantum groups 
on $C^*$--algebras, proved that, as in the classical case, actions
of compact quantum groups
on quotient spaces are ergodic. But he also showed that not all
ergodic actions are of this form. More precisely, he found an example of a 
compact quantum group acting ergodically  on a (even  commutative) 
$C^*$--algebra,
which is not a quotient action by a point stabilizer subgroup.

Wang's example suggests that the desired identification $X=K\backslash G$
should be relaxed, in the quantum case, to the possibility of finding a
faithful inclusion
of ${\cal C}(X)$ in  the quantum quotient space ${\cal C}(K\backslash G)$.

In the first sections we revisit the notion of compact  quantum subgroup and
quotient space of a compact quantum group. 
 In Sect. 2 we define  a compact quantum group $K=({\cal
A}',\Delta')$ to be a subgroup of 
 a compact quantum group $G=({\cal A},\Delta)$ if 
there exists 
a surjective $^*$--homomorphism $\pi:{\cal
A}\to{\cal A}'$ intertwining the coproducts. This definition
 reduces to
Podles's definition for
compact matrix pseudogroups and agrees with that considered by Wang \cite{Wang_free}, with slightly different terminology.
In the same paper,
Wang    introduced the notion 
of {\it Woronowicz $C^*$--ideal} and showed that such ideals     
characterize  closed ideals   of ${\cal A}$ that correspond to quantum subgroups.
We introduce the notion of a {\it closed coideal} of a Hopf $C^*$--algebra. This is a closed ideal which is also a coideal in the analytic sense,  in that the role of the   tensor product in the algebraic framework is replaced by the spatial tensor 
product.
Every closed coideal is a Woronowicz $C^*$--ideal, and it is not clear to us whether   the converse   holds in general (cf. Remark 3).
 We show that a Woronowicz $C^*$--ideal is a closed coideal    in the special case where the $C^*$--algebra describing the corresponding quantum subgroup   is  nuclear 
  (Theorem 2.8). The general theory of nuclear $C^*$--algebras (see, e.g. \cite{Blackadar}, \cite{Murphy}) guarantees that this is always the case if the Hopf $C^*$--algebra of $G$ itself is nuclear.  In the general case, looking at non coamenable quantum groups, we see that the smooth part ${\cal I}_\infty$ of a Woronowicz $C^*$--ideal ${\cal I}$  may not be dense in ${\cal I}$.
  However, ${\cal I}_\infty$ always turns out to be  an algebraic coideal (Lemma 3.7)
which  determines completely
the representation category of the corresponding 
subgroup $K$. 
 In particular,
  every Woronowicz $C^*$--ideal ${\cal I}$ contains a canonical closed coideal, the norm closure 
  $\bar{{\cal I}_\infty}$. We show that if a quantum subgroup   is coamenable then  
  the smooth
  part ${\cal I}_\infty$ of the associated Woronowicz $C^*$--ideal ${\cal I}$  is dense in it  (Cor. 3.8).
In the general case,   one can always replace a subgroup $K$ with another subgroup 
$K_{\text{max}}$
for which   the associated Woronowicz $C^*$--ideal  is a closed coideal with dense smooth part,  
and with the same representation category. We also show that
the equivalence classes of subgroups for which their smooth subalgebras
are isomorphic as Hopf $^*$--algebras
(equivalent subgroups), are in bijective correspondence with the Hilbert
space $C^*$--subcategories  with tensor products, subobjects and
direct sums, containing $\text{Rep}(G)$, the representation category of
$G$.

We next consider the notion of representation induced by a representation of a 
subgroup $K$ to the whole group $G$ and we show a Frobenius reciprocity
theorem at an algebraic level (Sect. 5).

In the next section we generalize Longo and Peligrad theorem to quantum
groups: we show that an action  of a compact quantum group
$G=({\cal
A},\Delta)$ on a 
unital $C^*$--algebra ${\cal C}$ is ergodic if and only if ${\cal C}$ has
no proper closed $G$--invariant left ideal ${\cal I}$ for which the left
ideal
generated by the image of ${\cal I}$ under the action is dense. 
In the case
where
the action and the Haar measure are faithful, we show that we can drop the
density assumption.
We give
more equivalent properties based on the lack of certain hereditary
$G$--invariant
$C^*$--subalgebras and certain open projections of ${\cal C}$.

In the last section we look for a necessary and sufficient condition
in order that an ergodic $G$--space ${\cal C}$ be embeddable in a quantum
quotient 
space by a  subgroup. The idea is the following. If such 
an embedding were possible, then, assuming for simplicity that the acting 
quantum group have an everywhere defined counit, the restriction of that
counit to the quotient space should induce a $^*$--character on ${\cal
C}$. Thus the existence of a $^*$--character is a necessary condition in
this case.
On the other hand it is not difficult to check that if $G$ is
a classical group, acting ergodically on ${\cal C}$, the existence of a
$^*$--character on ${\cal C}$ actually 
forces ${\cal C}$ to be commutative (see Prop. 7.1), and therefore 
a quotient $G$--space. 
In the quantum group case, assuming then the existence of a $^*$--character $\chi$ on ${\cal C}$,
we define the  quantum subgroup $G_\chi$ stabilizing $\chi$ 
and we
show that if 
the action is ergodic and faithful and if 
$G$ has faithful Haar measure,
then  ${\cal C}$ can be embedded faithfully,
with its
$G$--action, into the   quantum quotient space $G_\chi\backslash G$ (Theorems 7.3 and 7.4).

\end{section}

\begin{section} {Quantum  subgroups and closed  coideals}

We start  defining the notion of  compact quantum subgroup
of a compact quantum group.

Let $G=({\cal A}, \Delta)$ be a separable compact quantum group
in the sense of Woronowicz \cite{WLesHouches}.
Recall that this is a pair of a separable unital $C^*$--algebra
${\cal A}$ and a unital $^*$--algebra homomorphism
$\Delta: {\cal A}\to{\cal A}\otimes{\cal A}$ such that
\begin{description}
\item a) $\Delta\otimes\iota\circ\Delta=\iota\otimes\Delta\circ\Delta$,
with $\iota$ the identity map on ${\cal A}$,
\item b) the sets 
$\{b\otimes I\Delta(c), b,c\in{\cal A}\}$
$\{I\otimes b\Delta(c), b,c\in{\cal A}\}$
both span dense subspaces of ${\cal A}\otimes{\cal A}$.
\end{description}

The map $\Delta$ is usually called the coproduct of ${\cal A}$
and property a) is referred to as {\it coassociativity} of $\Delta$.

Let $G=({\cal A},\Delta)$ be a compact quantum group, and
$\theta({\cal A},{\cal A}):{\cal
A}\otimes{\cal A}\to{\cal A}\otimes{\cal A}$ the automorphism which
exchanges the order in the tensor product. Then
compact quantum group $G_{\text{o}}:=({\cal
A},\theta({\cal A},{\cal A})\circ\Delta)$ will be   referred to as
the
group
{\it opposite}
to $G$. 
\medskip

\noindent{\bf 2.1 Definition} 
A pair
$K=({\cal A}', \Delta')$ consituted by a unital $C^*$--algebra ${\cal A}'$ 
and a unital $^*$--homomorphism $\Delta':{\cal A}'\to{\cal A}'\otimes{\cal
A}'$
is a {\it compact quantum subgroup}
of $G=({\cal A},\Delta)$ if there exists a unital $^*$--epimorphism
$\pi: {\cal A}\to{\cal A}'$ such that
$\pi\otimes\pi\circ\Delta=\Delta'\circ\pi$.\medskip

\noindent{\bf 2.2 Proposition} {\sl A compact quantum subgroup of a
separable
compact quantum group is a separable compact quantum group.}\medskip

\noindent{\bf Proof}
 Since $\pi$ is surjective,   
${\cal A}'$ must be separable as well. Furthermore $\Delta'$ is
coassociative, as, for $a'=\pi(a)\in{\cal A}'$,
$$\Delta'\otimes\iota\circ\Delta'(a')=\Delta'\otimes\iota(\Delta'(\pi(a)))=$$
$$\Delta'\otimes\iota(\pi\otimes\pi(\Delta(a)))=
\Delta'\circ\pi\otimes\pi(\Delta(a))=$$
$$\pi\otimes\pi\otimes\pi(\Delta\otimes\iota\circ\Delta(a))=
\pi\otimes\pi\otimes\pi(\iota\otimes\Delta\circ\Delta(a))=
\iota\otimes\Delta'\circ\Delta'(a').$$
Furthemore 
the same intertwining 
relation between coproducts shows that 
property b) holds  for
$K$ as well. In fact, for example, for $b'=\pi(b),c'=\pi(c)\in {\cal A}'$,
elements of the form
$$I\otimes
b'\Delta'(c')=I\otimes\pi(b)\Delta'\circ\pi(c)=\pi\otimes\pi(I\otimes
b\Delta(c))$$ 
span a dense subspace of $\pi\otimes\pi({\cal A}\otimes{\cal A})={\cal
A}'\otimes{\cal A}'$.
Thus $K$ is a separable compact quantum group.
\medskip

One can easily recognize that in the case of compact matrix psudogroups,
this definition agrees with Podles definition \cite{Podles}.
\medskip

\noindent{\it Remark 1}
If the algebra ${\cal A}$ is commutative, then ${\cal A}'$ is commutative 
as well, since $\pi$ is surjective. Let $G$ and $K$ be the spectra of
${\cal A}$
and ${\cal A}'$ respectively,
which must be compact groups. The epimorphism $\pi$ then defines 
an injective continuous map $K\to G$, and therefore an identification 
of $K$ with a closed subgroup of $G$.
\medskip

\noindent{\it Remark 2}
It is known that the Haar measure $h$ of a compact quantum group
$G=({\cal A},\Delta)$
is faithful if and only if the corresponding GNS representation $\pi_h$ is
faithful \cite{Wcmp}. 
Consider the compact
quantum group 
$G_{h}=(\pi_h({\cal A}), \Delta_h)$ defined in \cite{Wcmp}.
The surjective map $\pi_h: {\cal A}\to\pi_h({\cal A})$ satisfies the
required
intertwining relation between $\Delta$ and $\Delta_h$.
Therefore in the case where $h$ is not faithful, according to the previous
definition, $G_{h}$ should be regarded as
a {\it proper} subgroup of $G$! 
\medskip

If $K=({\cal A}',\Delta')$ is a compact quantum subgroup of $G=({\cal
A},\Delta)$, the  $C^*$--algebra structure 
of ${\cal A}'$ is precisely the $C^*$--structure of the quotient
$C^*$--algebra ${\cal A}/\text{ker }\pi$. Thus the coproduct
of ${\cal A}'$ can be pulled back to a coproduct on ${\cal A}/\text{ker
}\pi$ making it into a compact quantum subgroup  of $G$ via the quotient
map.

Given a compact quantum group $G=({\cal A}, \Delta)$
what properties should a closed
ideal ${\cal I}$ of ${\cal A}$     satisfy in order that
${\cal A}/{\cal I}$ become a compact quantum subgroup of $G$ through
the quotient map $\pi:{\cal A}\to{\cal A}/{\cal I}$? 

In the case where ${\cal A}$ is just a Hopf algebra, the coproduct takes values in the algebraic tensor product ${\cal A}\odot{\cal A}$. In this case, it is easy to see, and we shall in fact see it later,  that the required condition on ${\cal I}$
 reduces precisely to the notion  of algebraic {\it coideal}, namely,
$$\Delta({\cal I})\subset {\cal I}\odot{\cal A}+{\cal A}\odot {\cal I}.$$
For compact quantum groups, one needs the notion of {\it Woronowicz $C^*$--ideal} introduced by    Wang.\medskip

\noindent{\bf 2.3 Definition} (\cite{Wang_free})
A Woronowicz $C^*$--ideal of $G=({\cal A}, \Delta)$ is a closed ideal ${\cal I}$ of ${\cal A}$
such that $$\Delta({\cal I})\subset\text{ker }(\pi\otimes\pi),$$ where $\pi: {\cal A}\to{\cal A}/{\cal I}$ is the quotient map.
\medskip

\noindent The answer to the above  question has then been given by  Wang.\medskip

\noindent{\bf 2.4 Theorem} (\cite{Wang_free}) {\sl Woronowicz  $C^*$--ideals correspond precisely to quantum subgroups.}\medskip 

In analogy with the algebraic case, the following notion is natural.\medskip

\noindent{\bf 2.5 Definition}
 Let ${\cal I}$ be a closed ideal of ${\cal A}$. We shall call
${\cal I}$ 
a {\it closed coideal} if
$$\Delta({\cal I})\subseteq
{\cal I}\otimes {\cal A}+{\cal A}\otimes{\cal I}.$$\medskip

We emphasize that $\otimes$ denotes the spatial tensor product between $C^*$--algebras. 
Recall 
that the sum of a
closed ideal and a
$C^*$--subalgebra in a $C^*$--algebra ${\cal B}$ is always a
$C^*$--subalgebra of ${\cal B}$
(Cor. 1.5.8 in \cite{Pedersen}). So the sum  
${\cal I}\otimes {\cal A}+{\cal A}\otimes{\cal I}$
is a closed ideal of ${\cal A}\otimes{\cal A}$.
\medskip

\noindent{\bf  Remark 3}
A closed coideal is clearly a Woronowicz $C^*$--ideal. Prop. 2.5 of the published version of this paper essentially claims 
equality of the two notions. However, the proof  there is unclear. 
I am deeply grateful to Shuzhou Wang for pointing this out to me. That unclarity originated from a 
gap of Lemma 2.4 of the published version, which we correct here below. 
As a result, up to date,  we show that the notions of a Woronowicz $C^*$--ideal and that of  a   closed coideal  coincide under additional assumptions.
\medskip

We give a proof of the following   known
fact (cf \cite{Murphy} Theorem 6.5.2) to clarify, via a reduction to an estimate, what one needs  
in order to show that
 a Woronowicz $C^*$--ideal is a closed coideal. 
If ${\cal A}$ is any $C^*$--algebra faithfully represented on a Hilbert
space ${\cal H}$, by $\iota_{\cal A}:{\cal A}\to{\cal
B}({\cal H})$
we shall denote the defining representation of  ${\cal A}$.
\medskip

\noindent{\bf 2.6 Lemma} {\sl Let ${\cal A}$ and ${\cal B}$ be
unital $C^*$--algebras, universally    represented on
Hilbert
spaces denoted by ${\cal H}$
and ${\cal H}'$ respectively, and let $\pi:{\cal A}\to{\cal B}({\cal K})$
be a 
Hilbert space representation of ${\cal A}$.

The kernel of  the  representation 
$$\pi\otimes\iota_{\cal B}: 
{\cal A}\otimes{\cal B}
\to{\cal B}({\cal
K}\otimes{\cal H}')$$
of the spatial tensor product
${\cal A}\otimes{\cal B}$, is $\text{ker }\pi\otimes{\cal B}$
if and only if
 $\text{ker }\pi$ admits a (positive, bounded) approximate identity $(u_\alpha)$ such that
for $a_1,\dots a_N\in{\cal A}$, $b_1,\dots, b_N\in{\cal B}$,
$$\lim_\alpha\|\sum_1^Na_i(I-u_\alpha)\otimes
b_i\|\leq  \|\sum_1^Na_i(I-P)\otimes b_i\|,$$ 
as operators on ${\cal H}\otimes{\cal H}'$, where $P$ denotes the strong limit of $(u_\alpha)$.

  In particular, this  holds   if  the algebraic tensor product
$\pi({\cal A})\odot{\cal B}$ has a unique $C^*$--norm (e.g. either $\pi({\cal A})$ or ${\cal B}$ is nuclear.)
}\medskip

\noindent{\bf Proof} Consider a positive bounded approximate identity
$\alpha\to u_\alpha$
of $\text{ker }\pi$. This is a net strongly convergent to the orthogonal 
projection $P$ on the closed subspace generated by 
$\text{ker }\pi{\cal H}$. Since that subspace is left invariant 
by ${\cal A}$, $P$ lies in the commutant of ${\cal A}$ in ${\cal B}({\cal
H})$. We thus have a $^*$--representation of $\pi({\cal A})$ on ${\cal H}$
defined by
$$\pi(a)\in\pi({\cal A})\to a(I-P)\in{\cal B}({\cal H}).$$
The projection $P$  is   the central projection of ${\cal A}''$ corresponding to the ultraweak
closure  $(\text{ker }\pi)''$ of $\text{ker }\pi$
via $(\text{ker }\pi)''={\cal A}''P$.
 If an element of $a\in{\cal A}$ satisfies
$$a(I-P)=0$$ then $$a\in(\text{ker }\pi)''\subset \text{ker }(\pi''), $$
hence $$\pi(a)=\pi''(a)=0.$$ Therefore
the representation 
$\pi(a) \to a(I-P)$
is  faithful, and hence isometric.
It follows that 
 for all $a\in{\cal A}$,
$$\lim_\alpha\|a-au_\alpha\|=\|\pi(a)\|=\|a(I-P)\|.$$
The tensor product representation 
$\pi\otimes\iota_{\cal B}$
from the spatial tensor product
${\cal
A}\otimes{\cal B}$ to operators on ${\cal K}\otimes{\cal H}'$, has
range $^*$--isomorphic to 
the spatial tensor product $\pi({\cal
A})\otimes{\cal B}$. 
Therefore, if $a_1,\dots, a_N\in{\cal A}$, $b_1,\dots,
b_N\in{\cal B}$, the quotient norm of $[\sum_1^Na_i\otimes b_i]$ with respect to the 
ideal $\text{ker }(\pi\otimes\iota)$ coincides with the norm of the image 
$\sum_1^N\pi(a_i)\otimes b_i$
in  $\pi({\cal
A})\otimes{\cal B}$.
By Turumaru's theorem (see, e.g. \cite{Sakai}, Prop. 1.22.9), asserting
that the tensor product
of faithful representations is faithful, and therefore isometric, the spatial norm 
of $\sum_1^N\pi(a_i)\otimes b_i$ can in turn be computed as
$$\|\sum_1^N\pi(a_i)\otimes b_i\|=\|\sum_1^Na_i(I-P)\otimes b_i\|,$$
where at the right hand side we are using the operator norm arising
from the Hilbert space ${\cal H}\otimes{\cal H}'$. Therefore
$$\|[\sum_1^Na_i\otimes b_i]\|=\|\sum_1^Na_i(I-P)\otimes b_i\|.$$
On the other hand,
$\alpha\to u_\alpha\otimes I$ is an approximate identity for
$\text{ker }\pi\otimes{\cal B}$ (spatial tensor product), hence  
the quotient norm of $[\sum_1^N a_i\otimes b_i]'$ with respect to the 
ideal $\text{ker }\pi\otimes{\cal B}$ equals
$$\|[\sum_1^N a_i\otimes b_i]'\|=
\lim_\alpha\|\sum_1^Na_i(I-u_\alpha)\otimes
b_i\|.$$
We note that we always have $$ \|\sum_1^Na_i(I-P)\otimes b_i\|\leq\lim_\alpha\|\sum_1^Na_i(I-u_\alpha)\otimes
b_i\|.$$ Taking into account the previous considerations, this inequality expresses boundedness of
    the  
$^*$--homomorphism of $C^*$--algebras
$${\cal A}\otimes{\cal B}/\text{ker }\pi\otimes{\cal B}\to{\cal
A}\otimes{\cal B}/\text{ker }(\pi\otimes\iota_{\cal B})$$
(i.e. the homomorphism naturally induced by the inclusion of closed ideals $\text{ker }\pi\otimes{\cal B}\subset\text{ker }(\pi\otimes\iota_{\cal B})$)
 on the dense $^*$--subalgebra ${\cal A}\odot{\cal B}/\text{ker }\pi\otimes{\cal B}$ of ${\cal A}\otimes{\cal B}/\text{ker }\pi\otimes{\cal B}$,
image of the algebraic tensor product  ${\cal A}\odot{\cal B}$ under the quotient map
${\cal A}\otimes{\cal B}\to{\cal A}\otimes{\cal B}/\text{ker }\pi\otimes{\cal B}$.
 Hence  the   requirement  in the statement,
i.e. the validity of the  reverse inequality,   is equivalent to the fact that  this homomorphism    is in fact isometric on that dense $^*$--subalgebra, and hence everywhere. This is thus also equivalent to the equality
$\text{ker }\pi\otimes{\cal B}=\text{ker }(\pi\otimes\iota_{\cal B})$ and 
the proof of the first part is complete.
 For the last part, we note that, since the natural inclusion of $^*$--algebras $\pi({\cal A})\odot{\cal B}\subset \pi({\cal A})\otimes{\cal B}$ is faithful, there are natural identifications
 $$({\cal A}\odot{\cal B})\cap\text{ker }(\pi\otimes\iota_{\cal B})=\text{ker }\pi\odot{\cal B}=
 ({\cal A}\odot{\cal B})\cap(\text{ker }\pi\otimes{\cal B}),$$ 
 hence the above
  $C^*$--homomorphism
   restricts to a $^*$--isomorphism between the dense $^*$--subalgebras
   $${\cal A}\odot{\cal B}/\text{ker }\pi\otimes{\cal B}\to{\cal
A}\odot{\cal B}/\text{ker }(\pi\otimes\iota_{\cal B})\simeq\pi({\cal A})\odot{\cal B}.$$
of 
${\cal A}\otimes{\cal B}/\text{ker }\pi\otimes{\cal B}$ and ${\cal
A}\otimes{\cal B}/\text{ker }(\pi\otimes\iota_{\cal B})$ respectively.
We thus have two $C^*$--norms on the algebraic tensor product $\pi({\cal A})\odot{\cal B}$
which must coincide if that $^*$--algebra has a unique $C^*$--norm.
\medskip

\noindent{\bf 2.7 Proposition} {\sl Let $K=({\cal A}',\Delta')$ be a 
compact quantum subgroup of $G=({\cal A},\Delta)$ defined by the
surjection $\pi: {\cal A}\to{\cal A}'$. If ${\cal A}'$ is a nuclear $C^*$--algebra (or, more generally, an exact $C^*$--algebra)  then the Woronowicz $C^*$--ideal $\text{ker }\pi$ is a closed ideal
and a closed coideal of ${\cal A}$.}\medskip

\noindent{\bf Proof}
For any $a\in\text{ker }\pi$, the element  $b:=\iota_{\cal
A}\otimes\pi(\Delta(a))\in{\cal A}\otimes{\cal A}'$
lies 
in the kernel of $\pi\otimes\iota_{{\cal A}'}$, as
$$\pi\otimes\iota_{{\cal
A}'}(b)=\pi\otimes\pi\Delta(a)=\Delta'(\pi(a))=0.$$
On the other hand by Lemma 2.6,
$$\text{ker }(\pi\otimes\iota_{{\cal A}'})=\text{ker }\pi\otimes{\cal
A}'=\text{Image} (\iota_{\text{ker }\pi}\otimes\pi),$$
since ${\cal A}'$ is assumed to be nuclear, hence
there is $c\in \text{ker }\pi\otimes{\cal A}$ such that
$$\iota_{\cal A}\otimes\pi(\Delta(a))=(\iota_{\text{ker
}\pi}\otimes\pi)(c).$$
Set $d:=\Delta(a)-c$. Since $$\iota_{\cal A}\otimes\pi(d)=(\iota_{\cal
A}\otimes\pi)(\Delta(a))-(\iota_{\text{ker }\pi}\otimes\pi)(c)=0,$$
$$d\in\text{ker }(\iota_{\cal A}\otimes\pi)={\cal A}\otimes\text{ker
}\pi,$$
where the last equality follows from Lemma 2.6 again.
So $$\Delta(a)\in{\cal A}\otimes\text{ker }\pi+\text{ker
}\pi\otimes{\cal
A}.$$
\medskip

\noindent{\it Remark} If ${\cal A}$ is a nuclear $C^*$--algebra then every quotient $C^*$-algebra 
of ${\cal A}$ is nuclear \cite{Murphy}, hence the notions of Woronowicz $C^*$--ideal and that of  closed coideal of $G=({\cal A}, \Delta)$ coincide.
\medskip

The following result provides an analogue of Wang's result, Theorem 2.4, for compact quantum groups with nuclear $C^*$--algebras.
\medskip

\noindent{\bf 2.8 Theorem} {\sl If ${\cal I}$ is a closed ideal and a
closed
coideal of ${\cal A}$, there is a  unique coproduct $\Delta^{\cal I}$ on
${\cal A}/{\cal I}$ making $K:=({\cal A}/{\cal I}, \Delta^{\cal I})$ into
a compact
quantum subgroup of $G$ through the quotient map $q:{\cal
A}\to{\cal A}/{\cal I}$. Any compact quantum subgroup $K=({\cal A}',\Delta')$ of $G$   for which   ${\cal A}'$ is nuclear,  is of this
form. }\medskip

\noindent{\bf Proof} The  the quotient $C^*$--algebra ${\cal A}':={\cal
A}/{\cal I}$
can be endowed with the unital $^*$--homomorphism
$$\widehat{\Delta}: {\cal A}'\to ({\cal A}\otimes{\cal A})/{\cal J},$$
with ${\cal J}:={\cal I}\otimes{\cal A}+{\cal A}\otimes{\cal I}$, 
taking the element $a+{\cal I}$ to the element $\Delta(a)+{\cal J}$.
On the other hand if $q:{\cal A}\to {\cal A}'$ is the quotient map,
then the $^*$--epimorphism $$q\otimes q: {\cal A}\otimes {\cal A}\to
{\cal
A}'\otimes{\cal
A}'$$
actually annihilates ${\cal J}$, and therefore it gives rise to a 
$^*$--epimorphism
$$\widehat{q\otimes q}:({\cal A}\otimes{\cal A})/{\cal J}\to{\cal
A}'\otimes{\cal
A}'.$$
Set 
$$\Delta^{\cal I}:=\widehat{q\otimes q}\circ\hat{\Delta}:{\cal A}'\to
{\cal
A}'\otimes
{\cal A}'.$$
One has: $$\Delta^{\cal I}(q(a))=q\otimes q(\Delta(a)).$$
By the previous proposition, any compact quantum subgroup of $G$ corresponding to a nuclear $C^*$--algebra is of
this form, with ${\cal I}$ the kernel of the defining surjective map.
\medskip

\end{section}

\begin{section} {Equivalent subgroups}

We have noticed (Remark 2) that the definition of quantum subgroup we
have given in the previous section has the disadvantage that the 
compact quantum group $G_{h}$ 
becomes a proper quantum subgroup of $G$ in the case where the Haar
measure $h$ is
not faithful. 
In this section we introduce the 
notion of equivalence between 
compact quantum subgroups which has the effect that $G_{h}$
becomes equivalent to the whole group $G$.

Recall that a $N$--dimensional unitary representation 
of a compact quantum group $G=({\cal A},\Delta)$ 
is a unitary matrix $u=(u_{i,j})\in \text{Mat}_N({\cal A})$
 such that $\Delta(u_{i,j})=\sum_1^Nu_{i,r}\otimes u_{r,j}$.
Let $H$ be an $N$--dimensional Hilbert space with orthonormal basis
$\{\psi_i,i=1,\dots,N\}$. 
Consider the right
Hilbert
${\cal A}$--module $H\otimes{\cal A}$, with inner product
$(\psi\otimes a, \phi\otimes b)_{{\cal A}}=(\psi,\phi)a^*b$, $a,b\in{\cal
A}$,
$\psi,\phi\in H$. 
The linear map $u:\psi_i\in H\to\sum_j\psi_j\otimes u_{j,i}\in
H\otimes{\cal A}$ satisfies:
$$(u(\psi),u(\phi))_{{\cal A}}=(\psi,\phi)_{{\mathbb C}}I, \psi,\phi\in
H,\eqno(3.1)$$
$$u\otimes\iota\circ u=\iota\otimes\Delta\circ u,\eqno(3.2)$$
$$u(H){\cal A}=H\otimes{\cal A}.\eqno(3.3)$$
Conversely, any linear map $u:H\to H\otimes{\cal A}$ satisfying
$(3.1)$--$(3.3)$ arises from a unitary representation of $G$.
 
If $H$ is infinite dimensional, the most general notion
of unitary representation has been given in \cite{BS}. More in
detail,
we consider the Hilbert space $L^2({\cal A})$ obtained completing
${\cal A}$ w.r.t. the inner product defined by the Haar measure.
A unitary representation of $G$ on $H$ is a unitary operator $W$ on the 
tensor product Hilbert space $H\otimes L^2({\cal A})$ satisfying
on $H\otimes L^2({\cal A})\otimes L^2({\cal A})$ the relation
$$W_{12}W_{13}V_{23}=V_{23}W_{12}.$$
Here $V$ denote the multiplicative unitary on $L^2({\cal A})\otimes
L^2({\cal A})$ defined by $V(a\otimes b)=\Delta(a)I\otimes b$.
Notice that  
equations $(3.1)$ and $(3.2)$
 make sense even in the case where $H$
is separable and infinite dimensional. We illustrate how one can get
a representation of $V$ from them. 

Let us consider
the natural continuous map $\tau$ from the right ${\cal A}$--Hilbert
module
$H\otimes{\cal A}$ to the Hilbert space $H\otimes L^2({\cal A})$.
A  map 
$u$ satisfying $(3.1)$ and  $(3.2)$
gives rise to a  map $W_u$ on
$H\otimes L^2({\cal A})$  defined by
$$W_{u}\psi\otimes a=\tau(u(\psi)a),\quad \psi\in H, a\in{\cal A}\subset
L^2({\cal A}). $$
We shall need the following fact in Sect. 5.\medskip

\noindent{\bf 3.1 Proposition} {\sl Let $H$ be a separable Hilbert
space, and $u:H\to H\otimes {\cal A}$ a linear
 map satisfying $(3.1)$ and $(3.2)$. If the set $\{u(\psi)a,
\psi\in H, a\in {\cal A}\}$ is total in $H\otimes{\cal A}$, the map
$W_u$ defined as above is a unitary representation of the 
multiplicative unitary $V$.}\medskip

\noindent{\bf 3.2 Definition} 
Let   $u$ be a unitary finite dimensional Hilbert space representation
of a compact quantum
group
$G$ on $H$, and let 
$K=({\cal A}',\Delta')$ be a compact subgroup of $G$, defined by
$\pi:{\cal A}\to{\cal A}'$. 
Then $u\upharpoonright_K:=\iota\otimes\pi\circ u:H\to H\otimes {\cal A}'$ 
is easily seen to be a unitary representation of $K$ on $H$, that
we call the {\it restriction} of $u$ to $K$.
\medskip

The space of {\it intertwining operators} $(u, u')$ between two unitary
representations on Hilbert spaces $H$ and $H'$ respectively, is the set of
all linear maps $T$ from $H$ to $H'$ such that
$$u'\circ T=T\otimes \iota \circ u,$$
with $\iota$ the identity map on ${\cal A}$.
The category $\text{Rep}(G)$ with objects finite dimensional unitary
$G$--representations and arrows
intertwining operators is known to be a tensor $C^*$--category with
conjugates \cite{Wcmp}.

Let us consider the functor from $\text{Rep}(G)$ to $\text{Rep}(K)$,
taking a representation $u$ of $G$ to the restricted representation 
$u\upharpoonright_K$, and acting trivially on the arrows. This is 
clearly a faithful tensor $^*$--functor.
If we identify each representation with its Hilbert space, 
this functor gives us an inclusion
$$\text{Rep}(G)\subset\text{Rep}(K)$$
of Hilbert space categories. \medskip

\noindent{\bf 3.3 Theorem} {\sl If $K$ is a compact quantum subgroup of
$G$,
the smallest full tensor $^*$--subcategory of $\text{Rep}(K)$
with
subobjects and direct sums containing all the restricted representations
$\{u\upharpoonright_K, u\in\text{Rep}(G)\}$,
is the whole $\text{Rep}(K)$.}\medskip

In order to prove this theorem we  shall need
 the smooth part  ${\cal A}_\infty$ of ${\cal A}$, the set of all linear
combinations of matrix
elements of all finite dimensional unitary representations of $G$. This is
a
dense $^*$--subalgebra of ${\cal A}$ such that $\Delta({\cal
A}_\infty)\subset
{\cal A}_\infty\odot{\cal A}_\infty$ (see \cite{WLesHouches}). Here
$\odot$
denotes the
algebraic tensor product.
\medskip

\noindent{\bf 3.4 Proposition} {\sl There are  choices of  complete sets
of 
inequivalent, unitary, irreducible  representations $\hat{G}=\{u^\alpha,
\alpha\in A\}$ and $\hat{K}=\{v^\beta,\beta\in B\}$ of $G$ and $K$
respectively,
such that for each $\alpha\in A$, $u^\alpha\upharpoonright_K$ splits 
into a  direct sum 
of elements of $\hat{K}$:
$u^\alpha\upharpoonright_K=m_1v^{\beta_1}\oplus\dots\oplus
m_Nv^{\beta_N}$, with
$\beta_1,\dots,\beta_N\in B$ and $m_1,\dots,m_N$
positive integers.
Every element of $\hat{K}$ arises in
this way. 
}\medskip

\noindent{\bf Proof}
Let us choose a complete set $\hat{G}=\{u^\alpha,\alpha\in A\}$ of
irreducible unitary 
representations of $G$. Since the representation coefficients
$u^\alpha_{i,j}$ span ${\cal A}_\infty$, which is dense in ${\cal A}$,
the set $\{\pi(u^\alpha_{i,j}),\alpha\in A,
i,j=1,\dots,\text{dim}(u^\alpha)\}$ span a dense $^*$--subalgebra of 
${\cal A}'$. 
On the other hand each $\pi(u^\alpha_{i,j})$ is the coefficient of the 
restricted
representation $u\upharpoonright_K$ of $u$ to  $K$. 
Let us consider a complete set $\hat{K}=\{v^\beta, \beta\in B\}$ of 
unitary irreducible representations of $K$.
Up to  replacing
$u^\alpha$ by a unitarily equivalent $G$--representation, 
we can assume that there exist $\beta_1,\dots, \beta_N\in B$ and
multiplicities $m_1,\dots, m_N\geq 1$ such that
$u^\alpha\upharpoonright_K=m_1v^{\beta_1}\oplus\dots\oplus
m_Nv^{\beta_N}$.
Since the set of all the coefficients of the $v$'s thus obtained is dense 
in ${\cal A}'$, when $\alpha$ ranges over $A$, we must obtain all the
irreducibles of $K$. 
\medskip

Theorem 3.3 is now an easy consequence of the previous proposition.
\medskip

\noindent{\bf 3.5 Definition} Let $G=({\cal A},\Delta)$ be a compact
quantum group. Two
compact quantum subgroups of $G$,
 $K_1=({\cal A}', \Delta')$ and $K_2=({\cal
A}'',\Delta'')$, defined by surjections $\pi_1:{\cal A}\to{\cal A}'$
and $\pi_2:{\cal A}\to{\cal A}''$ respectively,  will
be called {\it equivalent}
if  their smooth parts ${\cal A}'_\infty$ and ${\cal A}''_\infty$ are
isomorphic as Hopf $^*$--algebras.
\medskip

\noindent{\bf 3.6 Proposition} {\sl  $K_1$ and $K_2$ are equivalent
subgroups
of $G$ if and only if 
$\text{Rep}(K_1)=\text{Rep}(K_2)$ 
as Hilbert space categories.
}\medskip

\noindent{\bf Proof} It is easy to show that an isomorphism between the
smooth parts of two subgroups induces an isomorphism of tensor
$^*$--categories, between the
corresponding representation categories leaving fixed the
representation Hilbert spaces.
Conversely, if 
$\text{Rep}(K_1)=\text{Rep}(K_2)$
as Hilbert space categories, 
the smooth parts of $K_1$ and $K_2$ must be isomorphic as
Hopf $^*$--algebras,
as each of them is isomorphic to the 
Tannaka--Krein dual of that category,
by Woronowicz Tannaka--Krein duality
theorem \cite{Wtk}.
\medskip

Consider the smooth part of the kernel of $\pi:{\cal A}\to{\cal A}'$,
${\cal I}_\infty:=\ker\pi\cap{\cal A}_\infty$, which is clearly a
$^*$--ideal of ${\cal A}_\infty$. 
\medskip

\noindent{\it Remark 3} ${\cal I}_\infty$  {\it is not dense}
in $\text{ker }\pi$ in general. Indeed, consider the
case of the  quantum subgroup
$G_{h}$ 
defined by the  representation
$\pi_h:{\cal A}\to\pi_h({\cal A})$
of $G$. If $\pi_h$ is not faithful (or, equivalently, $h$ is not
faithful), 
${\cal I}_\infty=\{0\}$ by \cite{Wcmp}, but $\text{ker
}\pi_h\neq\{0\}$.
\medskip

\noindent{\bf 3.7 Lemma} {\sl ${\cal I}_\infty$ is a $^*$--ideal of ${\cal
A}_\infty$ and an algebraic coideal of $({\cal A}_\infty,\Delta)$, in the
sense
that 
$\Delta({\cal I}_\infty)\subset {\cal I}_\infty\odot{\cal
A}_\infty+{\cal A}_\infty\odot{\cal I}_\infty$.}\medskip

\noindent{\bf Proof}
The only nontrivial statement is that 
$\Delta({\cal I}_\infty)\subset {\cal I}_\infty\odot{\cal
A}_\infty+{\cal A}_\infty\odot{\cal I}_\infty$.
We choose $\hat{G}$ and $\hat{K}$ as in the previous theorem.
Let us call ``diagonal'' all the pairs $(i,j)$ corresponding 
to the entries in the matrix $u^\alpha$ 
corresponding to the entries occupied by the coefficients of the 
$v^{\beta_p}$'s in the matrix $u^\alpha\upharpoonright_K$,
and ``off diagonal'' 
the remaining pairs.
Let now $X\in{\cal I}_\infty$ be written uniquely as a linear
combination of the coefficients of a finite set $\{u^\alpha,\alpha\in F\}$
of
$G$--representations, i.e.
$$X=
\sum_{\alpha\in F}\sum_{i,j=1}^{\text{dim}
u^\alpha}\lambda^\alpha_{i,j} u^\alpha_{i,j}=$$
$$\sum_{\alpha\in F}\sum_{\text{diagonal pairs}}
\lambda^\alpha_{i,j} u^\alpha_{i,j}+
\sum_{\alpha\in F}\sum_{\text{off diagonal pairs}}
\lambda^\alpha_{i,j} u^\alpha_{i,j}=
$$
$$
X'+X''.
$$
Now, if $(i,j)$ is off diagonal,
$\pi(u^\alpha_{i,j})=0$, therefore
$u^\alpha_{i,j}\in\text{ker}\pi$, which shows that $X''$, 
and therefore also $X'$, lie in $\text{ker}\pi$
as well. Therefore we are left to show that both $\Delta(u^\alpha_{i,j})$,
for all  $(i,j)$ off diagonal, and $\Delta(X')$ belong to ${\cal
I}_\infty\odot{\cal A}_\infty+{\cal A}_\infty\odot{\cal I}_\infty$.
Since, for $(i,j)$ off diagonal, we can write
$$\Delta(u^\alpha_{i,j})=
\sum_{k:(i,k)\text{ is off diagonal}}u^\alpha_{i,k}\otimes u^\alpha_{k,j}+
\sum_{k: (i,k)\text{ is diagonal}} u^\alpha_{i,k}\otimes u^\alpha_{k,j},$$
we realize that $\Delta(u^\alpha_{i,j})\in{\cal I}_\infty\odot {\cal
A}_\infty+
{\cal A}_\infty\odot{\cal I}_\infty$ (notice that in the second sum 
each $(k,j)$ must be off diagonal, as $(i,j)$ is).

Let us now think of  the element $X'\in\text{ker}\pi$. 
Let $F'=\{\beta_1,\dots,\beta_N\}$  denote the finite subset of $B$
of all unitary irreducible $K$--representations obtained from the 
irreducible components of every
representation in the set $\{u^\alpha\upharpoonright_K,\alpha\in F\}$.
For each $r=1,\dots, N$ consider the subset $F_r\subset F$ constituted
by all $\alpha\in A$ for which $v^{\beta_r}$ is a subrepresentation of
$u^{\alpha}\upharpoonright_K$. For each such $\alpha$, let $m^\alpha_r$ 
be the multiplicity of $v^\beta_r$ in $u^\alpha\upharpoonright_K$.
We choose unitary equivalences so that for each $\alpha\in F_r$,
the subrepresentations of the form $m^\alpha_jv^{\beta_j}$ of 
$u^\alpha\upharpoonright_K$ are listed in order with increasing indices 
$j$. Then 
$$0=\pi(X')=\sum_{\alpha\in
F}\sum_{(i,j)\text{diagonal}}\lambda^\alpha_{i,j}\pi(u^\alpha_{i,j})=$$
$$\sum_{r=1}^N\sum_{h,k=1}^{\text{dim} v^\beta_r}
(\sum_{\alpha\in 
F_r}\sum_{p=1}^{m^\alpha_r}\lambda^\alpha_{i^\alpha_p,j^\alpha_p})
v^{\beta_r}_{h,k},$$
where, denoting by $f_j$ the dimension of $v^{\beta_j}$, and with 
$\delta^\alpha_r:=\sum_{j=1}^{r-1}m^\alpha_jf_j$, 
we have
set:
$i^\alpha_1=\delta^\alpha_r+h$, $i^\alpha_2=\delta^\alpha_r+f_r+h$,
$\dots$, $i^\alpha_{m^\alpha_r}=\delta^\alpha_r+(m^\alpha_r-1)f_r+h$,
$j^\alpha_1=\delta^\alpha_r+k$, $j^\alpha_2=\delta^\alpha_r+f_r+k$,
$\dots$, $j^\alpha_{m^\alpha_r}=\delta^\alpha_r+(m^\alpha_r-1)f_r+k$.
Since the $v^\beta_{h,s}$'s are linearly independent, we see that
for each possible $r$, and each possible pair $(h,k)\in\{1,\dots,
f_r\}\times\{1,\dots, f_r\}$,
$$\sum_{\alpha\in F_r}\sum_{p=1}^{m^\alpha_r}\lambda^\alpha_{i_p,j_p}=0.$$
Then 
$$X'=\sum_{r}\sum_{h,k} X^r_{h,k},$$
where
$$X^r_{h,k}=\sum_{\alpha\in
F_r}\sum_{p=1}^{m^\alpha_r}\lambda^\alpha_{i_p,j_p}u^\alpha_{i_p,j_p},$$
all of them belonging to $\text{ker}\pi$.
Therefore it suffices to show that $$\Delta(X^r_{h,k})\in{\cal
I}_\infty\odot{\cal A}_\infty+{\cal A}_\infty\odot{\cal I}_\infty$$
for each fixed $r, h, k$.
Set, for simplicity, $\mu^\alpha_p:=\lambda^{\alpha}_{i_p, j_p}$.
Then
$$\Delta(X^r_{h,k})=
\sum_{\alpha\in F_r}\sum_{p=1}^{m^\alpha_r}\sum_{q=1}^{\text{dim}
u^\alpha}\mu^\alpha_p
u^\alpha_{i_p, q}\otimes u^{\alpha}_{q,j_p}.$$
We split the sum in $q$ in two parts: the sum with $q$ ranging the
interval
$$\delta^\alpha_r+(p-1)f_r+1,\dots,
\delta^\alpha_r+(p-1)f_r+f_r=\delta^\alpha_r+pf_r$$
and the sum over the remaining values. We get
$$\sum_{\alpha\in
F_r}\sum_{p=1}^{m^\alpha_r}\sum_{q=\delta^\alpha_r+
(p-1)f_r+1}^{\delta^\alpha_r+(p-1)f_r+p}\mu^\alpha_p
u^\alpha_{i_p, q}\otimes u^{\alpha}_{q,j_p}+$$
$$\sum_{\alpha\in F_r}\sum_{p=1}^{m^\alpha_r}\sum_{\text{the remaining
indices } q  }\mu^\alpha_p
u^\alpha_{i_p, q}\otimes u^{\alpha}_{q,j_p}.$$
Each pair $(i_p, q)$ and $(q, j_p)$ arising from the second sum is off
diagonal, so the second sum belongs to ${\cal I}_\infty\odot{\cal
I}_\infty$.
Let us think of the first sum. We first perform the sum over the pairs
$(i_p,q)$ in the set
$$A_1:=\{(i_1, \delta^\alpha_r+1), (i_2,\delta^\alpha_r+f_r+1),\dots,
(i^{m^\alpha_r},\delta^\alpha_r+(m^\alpha_r-1)f_r+1)\},$$
followed by the sum over the pairs
$$A_2:=\{(i_1, \delta^\alpha_r+2), (i_2,\delta^\alpha_r+f_r+2),\dots,
(i^{m^\alpha_r},\delta^\alpha_r+(m^\alpha_r-1)f_r+2)\},$$
and, finally at the last step, over the pairs in the set
$$A_{f_r}:=\{(i_1, \delta^\alpha_r+f_r),
(i_2,\delta^\alpha_r+f_r+f_r),\dots,
(i^{m^\alpha_r},\delta^\alpha_r+(m^\alpha_r-1)f_r+f_r)\}.$$
We show then that each addendum
$$\sum_{\alpha\in F_r}\sum_{(i_p, q)\in
A_s}\mu^\alpha_pu^\alpha_{i_p,q}\otimes u^\alpha_{q,j_p}$$
lies in ${\cal I}_\infty\odot{\cal A}_\infty+{\cal A}_\infty\odot{\cal
I}_\infty$.
In fact, assuming for simplicity $s=1$, that sum can be also written as
$$(\sum_{\alpha\in F_r}\sum_{(i_p,q)\in A_1}\mu^\alpha_p
u^\alpha_{i_p,q})\otimes u^\alpha_{\delta^\alpha_r+1,j_1}+
\sum_{\alpha\in F_r}\sum_{(i_p,q)\in A_1}\mu^\alpha_p
u^\alpha_{i_p,q}\otimes
(u^\alpha_{q,j_p}-u^\alpha_{\delta^\alpha_r+1,j_1}),$$
and the claim is proved.
\medskip

\noindent{\bf 3.8 Corollary} {\sl A Woronowicz $C^*$--ideal ${\cal I}$ contains a canonical closed coideal, the norm closure of ${\cal I}_\infty={\cal I}\cap{\cal A}_\infty$.
If the associated subgroup $K=({\cal A}/{\cal I}, \Delta^{\cal I})$ is coamenable then
${\cal I}_\infty$ is dense in ${\cal I}$.
}
\medskip

\noindent{\bf Proof} The first statement follows from Lemma 3.7 and a routine  completeness argument. We next show that
of $\bar{{\cal I}_\infty}={\cal I}$.
Since $\bar{{\cal I}_\infty}$ and ${\cal I}$ have the same intersection with
${\cal A}_\infty$,
the quotient map ${\cal A}/\bar{{\cal I}_\infty}\to{\cal A}/{\cal I}$
restricts to a $^*$--isomorphism between the dense $^*$--subalgebras
 ${\cal A}_\infty/\bar{{\cal I}_\infty}\to{\cal A}_\infty/{\cal I}$. On the other hand,  
${\cal A}_\infty/{\cal I}$, can be identified with the canonical dense subalgebra of the quantum subgroup
${\cal A}/{\cal I}$, which has a unique $C^*$--norm by assumption. Hence
the two quotient norms need to coincide via the $^*$--isomorphism. Therefore
${\cal I}=\bar{{\cal I}_\infty}$.\medskip

\noindent{\it Remark 4}
Let us consider complete sets $\hat{G}$, $\hat{K}$ of irreducible
representations of $G$ and $K$ respectively, as in 
proposition 3.4.
We thus see that
 $\pi$ restricts to a surjective $^*$--homomorphism
${\cal
A}_\infty\to{\cal A}'_\infty$ with kernel ${\cal I}_\infty$.
Since ${\cal I}_\infty$ is an algebraic ideal and coideal of the Hopf
$^*$--algebra
$({\cal A}_\infty, \Delta)$, $({\cal A}'_\infty,\Delta')$ and $({\cal
A}_\infty/{\cal I}_\infty,\Delta^{{\cal I}_\infty})$ are isomorphic as
Hopf $^*$--algebras.
\medskip

\noindent{\it Remark 5} We show that up to replacing $K$ with an
equivalent subgroup, $K_{\text{max}}$, we can always assume that the
smooth part of the
${\cal I}_\infty$ be dense in $\text{ker }\pi$.
If $K=({\cal A}',\Delta')$ is a compact quantum subgroup of $G=({\cal
A},\Delta)$, we can complete the $^*$--ideal ${\cal I}_\infty=\text{ker
}\pi\cap{\cal A}_\infty$ and obtain a closed ideal 
 and a closed
coideal $\bar{{\cal I}_\infty}$.
Then thanks to Theorem 2.8 and Lemma 3.7, we can form another compact
quantum subgroup
$K_{\text{max}}$ by taking the quotient with respect to $\bar{{\cal 
I}_\infty}$.
Let $q$ be the corresponding quotient map.
This subgroup is clearly equivalent to the original subgroup $K$.
But $K$ can be in turn regarded as a subgroup of $K_{\text{max}}$,
as the inclusion $\bar{{\cal I}_\infty}\subset\text{ker }\pi$ provides
a 
$^*$--epimorphism $$\alpha: {\cal A}/\bar{{\cal I}_\infty}\to{\cal
A}/\text{ker
}\pi={\cal
A}'$$
such that, for $a\in{\cal A}$,
$$\Delta'\circ\alpha(a+\bar{{\cal
I}_\infty})=\Delta'(\pi(a))=\pi\otimes\pi\circ\Delta(a).$$
On the other hand
$$\alpha\otimes\alpha\circ\Delta^{\bar{{\cal I}_\infty}}(a+\bar{{\cal
I}_\infty})=\alpha\circ
q\otimes\alpha\circ q\Delta(a)=\pi\otimes\pi\circ\Delta(a),$$
so $\Delta'\circ\alpha=\alpha\otimes\alpha\circ\Delta^{\bar{{\cal
I}_\infty}}$.
\medskip

Let $h'$ be the Haar measure on $K$. The associated GNS representation
$\pi_{h'}:{\cal A}'\to\pi_{h'}({\cal A}')$ has as image the 
compact quantum group $K_{h'}$, which is a subgroup of
$K$. Notice that $\alpha$ is an isomorphism on 
 the smooth part of ${\cal A}/\bar{{\cal I}_\infty}$.
Therefore 
composing $h'\circ\alpha$
gives a positive state on
${\cal A}/\bar{{\cal I}_\infty}$ which acts as the Haar measure
$h'_{\text{max}}$ of
$K_{\text{max}}$. Thus
$L^2(K_{\text{max}},h'_{\text{max}})=L^2(K, h')$ and
$\pi_{h'}\circ\alpha=\pi_{h'_{\text{max}}}$.
These arguments show that if $h'_{\text{max}}$ is faithful then both 
$\pi_{h'}$ and $\alpha$ are faithful maps, so the equivalence class of
subgroups of $G$ equivalent to $K$ is constituted by all $^*$--isomorphic
quantum groups.
\medskip

We have seen  that if $K$ is a compact quantum subgroup of
$G$
then $\text{Rep}(G)\subset\text{Rep}(K)$ as a tensor $C^*$--category and
that equivalent subgroups have the same representation categories (Prop.
3.6).
Therefore there is a well defined map associating to the equivalence 
class $[K]$ of a subgroup $K$ of $G$ a  Hilbert space category, 
$\text{Rep}(K)$, containing $\text{Rep}(G)$.\medskip

\noindent{\bf 3.9 Theorem} {\sl Let $G$ be a compact quantum group. The
map
$[K]\to \text{Rep}(K)$ is a 
 bijective correspondence between equivalence classes of 
quantum subgroups  of $G$ and Hilbert space $^*$--categories with 
 tensor products,
 conjugates, direct sums and subobjects, containing
$\text{Rep}(G)$.}\medskip

\noindent{\bf Proof} We just need to show that the map is 
surjective. Let ${\cal T}$ be a tensor $^*$--category containing
$\text{Rep}(G)$ as in the statement. Then by Woronowicz Tannaka--Krein 
duality theorem we can find  a unital Hopf $^*$--algebra $K_0=({\cal
B},\Delta')$ with a unital
coassociative coproduct $\Delta'$ which is a universal model for ${\cal
T}$ (see \cite{Wtk}). On the other hand, since ${\cal
T}\supset\text{Rep}(G)$,
by universality of $({\cal A}_\infty,\Delta)$ there must exist a
surjective $^*$--homomorphism  ${\cal A}_\infty\to{\cal B}$ intertwining
the corresponding coproducts. The kernel ${\cal I}_0$  is a $^*$--ideal
and also a  
coideal (use the same arguments as in the proof of Lemma 3.7 to show that
it is a coideal). Complete ${\cal I}_0$ in ${\cal A}$ and obtain a closed
ideal and a
closed coideal ${\cal I}$ in ${\cal A}$. The corresponding quantum
subgroup $K$ of  
$G$ has then ${\cal T}$ as representation category.

 \end{section}

\begin{section} {Quantum quotient spaces}

Recall that an {\it action} of a compact quantum group $G=({\cal A},
\Delta)$ on a unital $C^*$--algebra ${\cal F}$ is a 
unital $^*$--homomorphism $\delta:{\cal F}\to{\cal F}\otimes{\cal A}$
such that
$$\delta\otimes\iota\circ\delta=\iota\otimes\Delta\circ\delta.$$
The map
$\delta$ will be called {\it ergodic} if the fixed point algebra
$${\cal F}^\delta:=\{f\in{\cal F}: \delta(f)=f\otimes I\}$$
reduces to the complex numbers.

We call a unital $^*$--homomorphism
$\delta':{\cal F}\to{\cal A}\otimes{\cal F}$ an {\it opposite action}
of $G$ if it satisfies the relation
$$\iota\otimes\delta'\circ\delta'=\Delta\otimes\iota\circ\delta'.$$
$\delta'$ is an opposite action of $G$ on ${\cal F}$ if and only if
$\delta:=\vartheta({\cal F},{\cal A})\circ\delta'$ is an action of the
opposite group $G_{\text{o}}$
on ${\cal F}$.
\medskip

\noindent{\bf 4.1 Proposition} {\sl
If $K=({\cal A}', \Delta')$ is a compact quantum subgroup of $G=({\cal A},
\Delta)$,
there is an action 
$$\delta_K:{\cal A}\to {\cal A}\otimes{\cal A}'$$ of $K$
on the $C^*$--algebra
${\cal A}$, defined by
$\delta_K:=\iota\otimes\pi\circ\Delta$.}\medskip

\noindent{\bf Proof} Indeed,
$$\delta_K\otimes\iota\circ\delta_K=
\iota\otimes\pi\otimes\iota\circ\Delta\otimes\iota\circ
\iota\otimes\pi\circ\Delta=$$
$$\iota\otimes\pi\otimes\pi\circ\Delta\otimes\iota
\circ\Delta=\iota\otimes\pi\otimes\pi\circ\iota\otimes\Delta\circ\Delta=$$
$$\iota\otimes\Delta'\circ\iota\otimes\pi\circ\Delta=
\iota\otimes\Delta'\circ\delta_K.$$\medskip

We then consider the fixed point algebra
$${\cal A}^{\delta_K}:=\{T\in {\cal A}: \delta_K(T)=T\otimes I\}$$
and call it {\it the quantum left coset space}.
Similarly, $\delta'_{K}:=\pi\otimes\iota\circ\Delta$
is an opposite action of $K$ on ${\cal A}$, whose fixed point
algebra
$${\cal A}^{\delta'_{K}}:=\{T\in{\cal A}: \delta'_K(T)=I\otimes T\}$$
will be called {\it the quantum right coset space}.
\medskip

Let $G$ be a compact group, and (${\cal C}(G)$, $\Delta$)  the associated
compact
quantum group
with coproduct $\Delta(f)(g,h)=f(gh)$. Let 
$K$ be a closed subgroup
of $G$. The  action of the  quantum subgroup group
$({\cal C}(K),\Delta')$ 
on ${\cal C}(G)$ just defined is given by 
$$\delta_K(f)(g,k)=f(gk), \quad g\in G,\quad k\in K.$$
Thus the fixed point algebra is 
$${\cal C}(G)^{\delta_K}=\{f\in {\cal C}(G): f(gk)=f(g), g\in G, k\in
K\},$$
the set of all class functions ${\cal C}(G/K)$. Its spectrum is the
compact
space $G/K$ of left cosets.
Similarly, $\delta'_K(f)(k, g)=f(kg)$, so 
$${\cal C}(G)^{\delta'_K}=\{f\in {\cal C}(G): f(kg)=f(g), g\in G, k\in
K\}={\cal C}(K\backslash G)$$
with spctrum the space $K\backslash G$ of right cosets.

These coset spaces are known to be endowed with left and right
$G$--actions by homeomorphisms. For example, for the right coset space,
$$g: Kg_1\in K\backslash G \to Kg_1 g\in K\backslash G$$
makes   $K\backslash G$  into a  homogeneous space,
in the sense that the above $G$--action is ergodic.
There is a natural way of defining the corresponding
$G$--actions in the quantum situation:
 restrict the coproduct   $\Delta$ to the 
left and right coset spaces. 
\medskip

\noindent{\bf 4.2  Proposition} {\sl {\rm \cite{Wang}} 
\begin{description}
\item{\rm a)}
The map 
$$\eta:=\Delta\upharpoonright_{{\cal A}^{\delta_K}}: {\cal
A}^{\delta_K}\to {\cal A}\otimes{\cal A}^{\delta_K}$$ 
is an ergodic opposite  action of $G$ on the quantum left coset space
${\cal
A}^{\delta_K}$ such that ${\cal A}\otimes I\eta({\cal A}^{\delta_K})$ 
is dense in ${\cal A}\otimes{\cal A}^{\delta_K}$.
\item{\rm b)}
The map $$\eta':=\Delta\upharpoonright_{{\cal A}^{\delta'_K}}:{\cal
A}^{\delta'_K}\to{\cal A}^{\delta'_K}\otimes{\cal A}$$
is an ergodic action of $G$ on the quantum right coset space ${\cal
A}^{\delta'_K}$ such that $I\otimes{\cal A}\eta'({\cal A}^{\delta'_K})$
is dense in ${\cal A}^{\delta'_K}\otimes{\cal A}$.
\end{description}}\medskip

\noindent{\bf Proof} We prove only a). Consider the map
$\iota\otimes\delta_K: {\cal A}\otimes{\cal A}\to {\cal A}\otimes{\cal
A}\otimes{\cal A}'$. This is an action of $K$ on ${\cal A}\otimes{\cal
A}$ with fixed point algebra ${\cal A}\otimes{\cal A}^{\delta_K}$. We show
that 
$\eta({\cal A}^{\delta_K})\subset{\cal A}\otimes{\cal A}^{\delta_K}$, or,
equivalently, that,
 for $T\in{\cal A}^{\delta_K}$,
 $\iota\otimes{\delta_K}(\Delta(T))=\Delta(T)\otimes I$.
Indeed,
$$\iota\otimes\delta_K(\Delta(T))=
\iota\otimes\iota\otimes\pi\circ\iota\otimes\Delta\circ\Delta(T)=$$
$$\iota\otimes\iota\otimes\pi\circ\Delta\otimes\iota
\circ\Delta(T)=
\Delta\otimes\iota\circ\delta_K(T)=$$
$$
\Delta\otimes\iota(T\otimes I)=
\Delta(T)\otimes I.$$
We  show that
the $\eta$--fixed point algebra reduces to the complex numbers.
Let $T\in{\cal A}^{\delta_K}$ satisfy
$\eta(T)=I\otimes T$. Since $\delta_K(T)=T\otimes I$, we  have both:
$$\Delta(T)=I\otimes T, \iota\otimes\pi\circ\Delta(T)=T\otimes I,$$
so
$$I\otimes\pi(T)=T\otimes I,$$ which shows that $T\in{\mathbb C}$.
We finally show that the linear span of elements of the form $a\otimes
I\eta(b')$, with $a\in{\cal A}$, $b'\in{\cal A}^{\delta_K}$, is dense.
Let $h'$ be the Haar measure of  $K$, and $E:{\cal A}\to{\cal
A}^{\delta_K}$ the conditional expectation: $E(b):=\iota_{\cal A}\otimes
h'\circ\delta_K$. A straightforward computation shows that $E$ commutes
with the action $\eta$, in the sense that
$$\eta\circ E=\iota_{\cal A}\otimes E\circ\Delta.$$
Take $b'\in{\cal A}^{\delta_K}$ of the form $E(b)$, with $b$ ranging
over ${\cal A}$. Then 
$$a\otimes I\eta(b')=\iota_{\cal A}\otimes E(a\otimes I\Delta(b)),$$
and the conclusion follows from the fact that ${\cal
A}\otimes I\Delta({\cal A})$ is dense in ${\cal A}\otimes{\cal A}$. 
\medskip

The  pairs $$G/K:=({\cal A}^{\delta_K}, \eta)$$ and 
$$K\backslash G:=({\cal A}^{\delta'_K}, \eta')$$
will be called {\it the compact quantum left and right 
quotient
spaces}, respectively,
defined by  the subgroup $K$ of $G$.

\end{section}
\begin{section}{ Induced  representations and Frobenius reciprocity}

In this section we 
define the representation of a compact quantum group induced
by a representation  of a compact quantum subgroup.

 Let 
$K=({\cal A}',\Delta')$ be a compact subgroup of $G$, defined by
$\pi:{\cal A}\to{\cal A}'$. 
and 
let  $u: H\to H\otimes{\cal A}'$ be a finite dimensional unitary
representation
of $K$ on the Hilbert space $H$.
Consider the opposite action of $K$ on ${\cal A}$: 
$$\delta'_{K}:=
\pi\otimes\iota\circ\Delta: {\cal A}\to{\cal
A}'\otimes{\cal A}.$$ 
Ona has:
$\delta'_K({\cal A}_\infty)\subset \pi({\cal A}_\infty)\odot{\cal
A}_\infty\subset{\cal A}'_\infty\odot{\cal A}_\infty$.
Define maps:
$$u\otimes\iota: H\otimes{\cal A}\to H\otimes{\cal A}'\otimes {\cal A},$$
$$\iota\otimes\delta'_{K}: H\otimes{\cal A}\to H\otimes{\cal
A}'\otimes{\cal A}.$$
Note that $u\otimes\iota(H\odot{\cal A}_\infty)\subset H\odot {\cal
A}'_\infty\odot{\cal A}_\infty$ and 
$\iota\otimes\delta'_K(H\odot{\cal A}_\infty)\subset H\odot
{\cal A}'_\infty\odot{\cal A}_\infty$.
Consider the following subspaces:

$$X_{\text{Ind}_\infty(u)}:= 
\{\xi\in H\odot{\cal A}_\infty:
u\otimes\iota(\xi)=\iota\otimes\delta'_K(\xi)\},$$
$$X_{\text{Ind}_0(u)}:= 
\{\xi\in H\otimes{\cal A}:
u\otimes\iota(\xi)=\iota\otimes\delta'_K(\xi)\},$$
the latter being a Banach subspace of $H\otimes{\cal A}$.
Notice that if $\xi,\eta\in X_{\text{Ind}_0(u)}$,
$(\xi,\eta)_{\cal A}$ is an element of the right coset space
${\cal A}^{\delta'_K}$: 
$$\delta'_{K}((\xi,\eta)_{\cal
A})=(\iota_H\otimes\delta'_K(\xi),\iota_H\otimes\delta'_K(\eta))_{{\cal
A}'\otimes{\cal A}}=$$
$$(u\otimes \iota_{\cal A}(\xi),u\otimes\iota_{\cal
A}(\eta))_{{\cal A}'\otimes{\cal A}}=I\otimes(\xi,\eta)_{\cal A}$$
by $(3.1)$.
Therefore $X_{\text{Ind}_0(u)}$ is a right Hilbert ${\cal
A}^{\delta'_K}$--module.

We shall also consider the norm closure $H_{\text{Ind}(u)}$ in $H\otimes
L^2({\cal A})$
of the image of $X_{\text{Ind}_\infty(u)}$ under the natural
continuous
map
$\tau:
H\otimes{\cal A}\to H\otimes
L^2({\cal A})$.

One should point
out that there do exist  functions in $H_{\text{Ind}(u)}$.
Indeed, let $$E:H\odot{\cal A}_\infty\to H\odot {\cal
A}_\infty$$ denote the linear map that takes 
a simple tensor $x=\psi\otimes a$ to the
element
$$E(x):=\iota_H\otimes (h'\circ m)\otimes\iota_{{\cal
A}}((\iota_H\otimes\kappa'\circ u)\otimes\delta'_K(x))=$$
$$\iota_H\otimes (h'\circ m)\otimes\iota_{{\cal A}}
((\iota_H\otimes\kappa'\circ
u)(\psi)\otimes\delta'_K(a)),
$$
where 
$m:{\cal A}'_\infty\odot{\cal A}'_\infty\to{\cal A}'_\infty$ denotes
the multiplication map: $a\otimes b\to ab$, $h'$ the Haar measure for
$K$ and $\kappa'$ the coinverse for ${\cal A}'_\infty$.
\medskip

\noindent{\bf 5.1 Lemma} {\sl 
$E$ is an idempotent map  with  $E(H\odot{\cal
A}_\infty)=
X_{\text{Ind}_\infty(u)}$.}\medskip

\noindent{\bf Proof} We show that, for any 
$x\in
X_{\text{Ind}_\infty(u)}$, $E(x)=x$.
$$E(x)=
\iota_H\otimes (h'\circ m)\otimes\iota_{{\cal
A}}((\iota_H\otimes\kappa'\circ u)\otimes\iota_{{\cal
A}'}\otimes\iota_{{\cal A}}\circ u\otimes\iota_{{\cal A}}(x)=$$
$$\iota_H\otimes (h'\circ m)\otimes\iota_{{\cal
A}}(\iota_H\otimes\kappa'\otimes\iota_{{\cal
A}'}\otimes\iota_{{\cal
A}}\circ\iota_H\otimes\Delta'\otimes\iota_{{\cal A}}\circ
u\otimes\iota_{{\cal A}}(x)=$$
$$\iota_H\otimes e'\otimes \iota_{{\cal A}}\circ u\otimes\iota_{{\cal
A}}(x)= x.$$
since $m\circ\kappa'\otimes\iota_{{\cal A}'}\circ\Delta'=e'$, with $e'$
the
counit of ${\cal A}'_\infty$, and $\iota_H\otimes e'\circ u=\iota_H$.
We are left to show that for any $x\in H\odot{{\cal A}_\infty}$,
$E(x)\in
X_{\text{Ind}_\infty(u)}$.
$$\iota_H\otimes\delta'_K(E(x))=
\iota_H\otimes (h'\circ m)\otimes\iota_{{\cal
A}'}\otimes\iota_{{\cal
A}}\circ\iota_H\otimes\iota_{{\cal
A}'}\otimes\iota_{{\cal A}'}\otimes\delta'_K((\iota_H\otimes\kappa'\circ
u)\otimes\delta'_K(x))=$$
$$\iota_H\otimes (h'\circ m)\otimes\iota_{{\cal
A}'}\otimes\iota_{{\cal
A}}((\iota_H\otimes\kappa'\circ
u)\otimes(\iota_{{\cal A}'}\otimes\delta'_K\circ\delta'_K)(x))=$$
$$\iota_H\otimes (h'\circ m)\otimes\iota_{{\cal
A}'}\otimes\iota_{{\cal
A}}((\iota_H\otimes\kappa'\circ
u)\otimes(\Delta'\otimes\iota_{{\cal A}}\circ\delta'_K)(x)).$$
Now for $x=\psi_i\otimes a$ with $a\in{\cal A}_\infty$ and
$\psi_i$ an element of an orthonormal basis of $H$, write 
$\delta'_K(a)=\sum a'_{1}\otimes a_{2}$ and $u(\psi_i)=
\sum_r \psi_r\otimes u_{r,i}$. Inserting these computations in the last
term gives, using strong right invariance of the Haar measure \cite{Wcmp}:
$\kappa'(h'\otimes\iota(a\otimes
I\Delta'(b)))=h'\otimes\iota(\Delta'(a)b\otimes I)$,
$$\sum_r\psi_r\otimes (h'\otimes\iota_{{\cal A}'}(\kappa'(u_{r,i})\otimes
I\Delta'(a'_{1})))\otimes a_{2}=$$
$$\sum_r\psi_r\otimes ({\kappa'}^{-1}(h'\otimes\iota_{{\cal
A}'}(\Delta'\circ\kappa' (u_{r,i})a'_{1}\otimes I)))\otimes a_{2}=$$
$$\sum_r\psi_r\otimes({\kappa'}^{-1}(\iota\otimes
h'(\kappa'\otimes\kappa'\circ\Delta'(u_{r,i})I\otimes a'_1)))\otimes
a_2=$$
$$\sum_r\psi_r\otimes(\iota\otimes
h'(\iota\otimes\kappa'\circ\Delta'(u_{r,i})I\otimes a'_1))\otimes a_2=$$
$$\sum_{r,k}\psi_r\otimes u_{r,k}\otimes h'(\kappa'(u_{k,i})a'_1)
a_2=$$
$$\sum_k u(\psi_k)\otimes h'(\kappa'(u_{k,i}a'_1))a_2=u\otimes
\iota_{{\cal
A}}(E(x)),$$
since 
$\Delta'\circ\kappa'=\vartheta\circ\kappa'\otimes\kappa'\circ\Delta'$
with $\vartheta$ the flip automorphism of ${\cal A}'\otimes{\cal A}'$.
\medskip

\noindent{\bf 5.2  Lemma} {\sl Let 
$(f_z)_{z\in{\mathbb C}}$ and
$(f'_z)_{z\in{\mathbb Z}}$ be the family of linear multiplicative 
functionals defined on ${\cal A}_\infty$ and ${\cal A}'_\infty$
respectively \cite{Wcmp}, and
let $u$ and $v$ be irreducible unitary representations of $K$ and $G$
respectively. If $v\upharpoonright_K=\oplus_1^m u\oplus u'$ 
with 
 $u'$
disjoint from $u$, 
then 
for all $r,s=1,\dots, N_u$, with $N_u$ the dimension of $u$, 
and $j=0,\dots, m-1$, there exist positive constants $\lambda_j$ such that
$$f_{-1}(v_{r+jN_u,s+jN_u})=\lambda_jf'_{-1}(u_{r,s}).$$}\medskip

\noindent{\bf Proof}  Let $F_v$ be the unique positive intertwiner
from $v$ to the double contragradient representation $v^{cc}$ with 
$\text{Tr}(F_v)=\text{Tr}(F_v^{-1})$. It is easy to check that
$F_v$ is an intertwiner from $v\upharpoonright_K$ to
$v\upharpoonright_K^{cc}$. Therefore $F_v$ leaves globally invariant
the space of $\oplus_1^m u$ and $F_j=E_jF_v\upharpoonright_{H_j}$ are
positive invertible
intertwiners from 
$u$ to $u^{cc}$, with $E_j$ the orthogonal projection from $H_v$
to the $j$--th space $H_j$ of $u$.
Therefore for some positive constants $\lambda_j$, $\lambda_jF_j=
F_u$.
By definition of $f_{-1}$ (see  (5.22) in \cite{Wcmp}), following
Woronowicz
notation, 
$$f_{-1}(v_{r+jN_u,
s+jN_u})=
\text{Tr}(F_v^{-1} m^v_{r+jN_u, s+jN_u})=$$
$$\text{Tr}(F_j^{-1} m^u_{r, s})=
\lambda_j\text{Tr}(F_u^{-1} m^u_{r,
s})=\lambda_jf'_{-1}(u_{r,s}).$$
\medskip

We next show that $E$ is a selfadjoint projection.\medskip

\noindent{\bf 5.3 Proposition} {\sl Let us regard $H\odot{\cal
A}_\infty$ as a dense subspace of the Hilbert space $H\otimes L^2({\cal
A})$. Then the densely defined operator $E: H\odot{\cal A}_\infty\to
H\otimes L^2({\cal A})$ is Hermitian:
$$(x', E(x))=(E(x'), x),\quad x,x'\in H\odot{\cal A}_\infty.\eqno(5.1)$$
Therefore $E$ is bounded and  extends uniquely to 
the orthogonal projection of $H\otimes L^2({\cal A})$ onto
$H_{\text{Ind}(u)}$.
}\medskip

\noindent{\bf Proof}
${\cal A}_\infty$ is linearly spanned by coefficients of unitary
irreducible representations of $G$ \cite{Wcmp}, therefore it suffices to
take
$x=\psi_i\otimes v_{p,q}$ and $x'=\psi_{i'}\otimes v'_{p',q'}$,
with $v$ and $v'$ irreducible representations of $G$ and $\{\psi_r\}$ an
orthonormal basis of $H$. 
A straightforward computation shows 
that for $x=\psi_i\otimes v_{p,q}$,
$$E(x)=\sum_{r,k}\psi_r\otimes h'(\kappa'(u_{r,i})\pi(v_{p,k}))v_{k,q}=
\sum_{r,k}h'(u_{i,r}^*\pi(v_{p,k}))\psi_r\otimes v_{k,q}.\eqno(5.2)$$
Since any unitary representation of $K$ is the direct sum
of irreducibles and since $E$ leaves globally invariant any subspace of
the 
form $H'\odot{\cal A}_\infty$, with $H'$ the space of a
subrepresentation of $u$, it suffices to assume $u$ irreducible.  
By the Peter-Weyl theory for compact quantum groups \cite{Wcmp},
$h'(u_{i,r}^*\pi(v_{p,k}))=0$ unless $v\upharpoonright_K$ contains
$u$ as a subrepresentation. 
In that case the computation of $E(x)$ shows that $(x', E(x))=0$ unless
$v'\upharpoonright_K$ contains $u$,
again by the Peter-Weyl theory of $K$.
Therefore both sides of $(5.1)$ annihilate, and therefore coincide,
unless both $v$ and $v'$ contain $u$ as a subrepresentation when
restricted to $K$. Let us assume then that this is the case.
The computation of $E(x)$ shows that, up to replacing 
$v$ by an equivalent representation, we can assume that
$v\upharpoonright_K$ takes the form:
$v\upharpoonright_K=\oplus_1^m u\oplus u'$, with $m$
the multiplicity of $u$ in $v\upharpoonright_K$.
If  $p$  is bigger than $mN_u$, with $N_u$ the dimension of
$u$, $E(x)=0$ since $u'$ and $u$ are disjoint.
For $p=jN_u+1,\dots,(j+1)N_u$, 
 for some $j=0,\dots m-1$,
$\pi(v_{p,k})=0$ for all $k$, unless
$k=jN_u+1,\dots,(j+1)N_u$. In this case $\pi(v_{p,k})=u_{p-jN_u, k-jN_u}$.
Then by \cite{Wcmp}, Theorem 5.7,
$$h(u_{i,r}^*\pi(v_{p,k}))=h(u_{i,r}^*u_{p-jN_u,k-jN_u})=\delta_{r,k-jN_u}
\frac{f'_{-1}(u_{p-jN_u,i})}{f'_1(\chi_u)},$$
with $z\to f_z$ the linear multiplicative functionals of ${\cal
A}'_\infty$ defined in \cite{Wcmp}, Theorem 5.6.
Therefore 
$$E(x)=
\frac{f'_{-1}(u_{p-jN_u,i})}{f'_1(\chi_u)}
\sum_r\psi_r\otimes
v_{r+jN_u,q}.$$
Here $\chi_u=\sum_s u_{s,s}$ is the character of the representation $u$.
Thus $E(x)=0$ unless $q=jN_u+1,\dots,(j+1)N_u$. Assume then that this is
the case.
Now 
$$(x',E(x))=\frac{f'_{-1}(u_{p-jN_u,i})}{f'_1(\chi_u)}\sum_r\delta_{i',r}
h({v'_{p',q'}}^*v_{r+jN_u,q})=$$
$$\frac{f'_{-1}(u_{p-jN_u,i})}{f'_1(\chi_u)}h({v'_{p',q'}}^*v_{i'+jN_u,q})=$$
$$\delta_{v',v}\delta_{q',q}\frac{f'_{-1}(u_{p-jN_u,i})}{f'_1(\chi_u)}
\frac{f_{-1}(v_{i'+jN_u,p'})}{f_1(\chi_v)}.$$
Exchanging the roles of $x$ and $x'$
gives
$$(E(x'),x)=\overline{(x,E(x')}=$$
$$\delta_{v,v'}\delta_{q,q'}\frac{\overline{f'_{-1}(u_{p'-j'N_u,i'})}}{\overline{f'_1(\chi_u)}}
\frac{\overline{f_{-1}(v_{i+j'N_u,p})}}{\overline{f_1(\chi_{v'})}}=$$
$$\delta_{v,v'}\delta_{q,q'}\frac{f'_{-1}(u_{i',p'-j'N_u,})}{f'_1(\chi_u)}
\frac{f_{-1}(v_{p,i+j'N_u})}{f_1(\chi_{v'})},$$
by the relation $\overline{f_{-1}(u_{r,s})}=f_{-1}(u_{s,r})$
shown in \cite{Wcmp}, $(5.17)$--$(5.18)$.
It suffices to assume then $v=v'$ and $q=q'$, so we also have $j=j'$.
The previous lemma now completes the
 proof of $(5.1)$.

Finally, $(5.1)$ shows that
$x-E(x)$ and $E(x)$ are orthogonal to each other
since $E^2=E$ by Lemma 5.1. Therefore 
in the Hilbert space norm: $\|E(x)\|^2\leq\|x\|^2$.
The rest is now clear.
\medskip

A  representation $v$ of $G$ on a vector space $V$ is  a 
linear map $v: V\to V\odot{\cal A}_\infty$ such that 
$$v\otimes\iota_{{\cal
A}}\circ v=\iota_V\otimes\Delta\circ v,\eqno(5.3)$$
$$v(V){\cal A}_\infty=V\odot{\cal A}_\infty.\eqno(5.4)$$
In the last equation $V\odot{\cal A}_\infty$ is regarded as a right
${\cal A}_\infty$--module in the natural way.
Note that a unitary finite dimensional representation of $G$ 
 is a vector space
representation of $G$.
Moreover, for any finite dimensional Hilbert space $H$,
$$\iota_H\otimes\Delta: H\odot{\cal A}_\infty\to H\odot{\cal
A}_\infty\odot{\cal A}_\infty$$ is vector space representation 
of $G$ on $H\odot{\cal A}_\infty$.

If $v$ and $v'$ are two vector space representations
of $G$
on  $V$ and $V'$
respectively, a  linear
map
$T:V\to V'$ is called an intertwiner 
if 
$$T\otimes 1_{{\cal A}}\circ v=v'\circ T
.$$
The space of all such intertwiners will be denoted by 
$(v,v')$.

We are now ready to define induced $G$--representations on the spaces 
$X_{\text{Ind}_\infty(u)}$ and $H_{\text{Ind}(u)}$.
\medskip

\noindent{\bf 5.4 Proposition} {\sl 
Let $u$ be a unitary finite dimensional representation of a compact
quantum subgroup $K$ of a compact quantum group $G$.
\begin{description}

\item {\rm a)} 
Then the idempotent $E: H\otimes{\cal A}_\infty\to
X_{\text{Ind}_\infty(u)}$ intertwines the vector space $G$--representation
$$\iota_H\otimes\Delta: H\odot{\cal A}_\infty \to H\odot{\cal
A}_\infty\odot{\cal
A}_\infty$$ 
with itself.
Therefore the restriction of $\iota_H\otimes \Delta$ to
$X_{\text{Ind}_\infty(u)}$ gives rise to a 
vector space representation
 $$\text{Ind}_\infty(u):
X_{\text{Ind}_\infty(u)}\to
X_{\text{Ind}_\infty(u)}\odot{\cal A}_\infty.$$

\item{\rm b)} The same map $\iota_H\otimes\Delta$ restricted to 
the right Hilbert module $X_{\text{Ind}_0(u)}$ 
is a bounded linear map $$\text{Ind}_0(u): X_{\text{Ind}_0(u)}\to
X_{\text{Ind}_0(u)}\otimes {\cal A},$$
where $\otimes$ denotes the exterior tensor product.

\item{\rm c)} 
The map $\text{Ind}_0(u)$  extends uniquely to a map
$$\text{Ind}(u): H_{\text{Ind}(u)}\to
H_{\text{Ind}(u)}\otimes{\cal A}$$
satisfying the assumptions of Prop. 3.1, and therefore a
unitary $G$--representation.
\end{description}
}\medskip

\noindent{\bf Proof} a) 
$H\odot{\cal A}_\infty$ is linearly spanned by simple tensors 
of the form $x=\psi_i\otimes v_{p,q}$, with $\{\psi_r\}$ an orthonormal
basis of $H$ and $v_{p,q}$ coefficients of an irreducible unitary
representation of $G$. By $(5.2)$
$$\iota_H\otimes\Delta\circ E(x)=
\sum_{r,k}
h'(u_{i,r}^*\pi(v_{p,k}))
\psi_r\otimes\Delta(v_{k,q})=$$
$$\sum_{r,k,s} h'(u_{i,r}^*\pi(v_{p,k}))
\psi_r\otimes
v_{k,s}\otimes v_{s,q}=
\sum_{s}E\otimes\iota(\psi_{i}\otimes v_{p,s}\otimes v_{s,q})=$$
$$E\otimes\iota\circ\Delta(x).$$
Therefore the restriction $\text{Ind}_\infty(u)$ of $\iota_H\otimes\Delta
$ to 
$X_{\text{Ind}_\infty(u)}$ 
takes that subspace into $X_{\text{Ind}_\infty(u)}\odot{\cal A}_\infty$ 
and clearly satisfies relation $(5.3)$.
We show $(5.4)$ with $V=X_{\text{Ind}_\infty(u)}=E(H\otimes{\cal
A}_\infty)$.
$$\text{Ind}_\infty(u)(E(H\otimes{\cal A}_\infty)){\cal A}_\infty=
E\otimes\iota(\iota_H\otimes\Delta(H\otimes{\cal A}_\infty)){\cal
A}_\infty=$$
$$E\otimes\iota(H\otimes(\Delta({\cal A}_\infty)I\otimes{\cal
A}_\infty)))=
E\otimes\iota(H\otimes{\cal A}_\infty\odot{\cal A}_\infty)=$$
$$X_{\text{Ind}_\infty(u)}\odot{\cal A}_\infty.$$

b) Obviously $\iota_H\otimes\Delta$ is a bounded linear map from the right
Hilbert module
$H\otimes{\cal A}$ to the exterior tensor product of right
Hilbert modules $(H\otimes{\cal
A})\otimes{\cal A}$. The space $X_{\text{Ind}_0(u)}$ is a norm closed
subspace of $H\otimes{\cal A}$. The following computations show 
that
$\iota_H\otimes\Delta(X_{\text{Ind}_0(u)})\subset
X_{\text{Ind}_0(u)}\otimes{\cal A}$. 
We  show that for  any $T\in
X_{\text{Ind}_0(u)}$,
$$u\otimes\iota_{\cal A}\otimes\iota_{\cal
A}(\iota_{H}\otimes\Delta(T)))=
\iota_H\otimes \delta'_{K}\otimes\iota_{\cal
A}(\iota_{H}\otimes\Delta(T))).
$$
The left hand side equals
$$\iota_H\otimes\iota_{\cal A'}\otimes\Delta(u\otimes\iota_{\cal A}(T))=
\iota_H\otimes\iota_{\cal
A'}\otimes\Delta(\iota_H\otimes \delta'_K(T))=$$
$$\iota_H\otimes\iota_{\cal
A'}\otimes\Delta\circ\iota_H\otimes\pi\otimes\iota_{\cal
A}\circ\iota_H\otimes \Delta(T)=
\iota_H\otimes\pi\otimes\iota_{\cal A}\otimes\iota_{\cal
A}\circ\iota_H\otimes\Delta\otimes\iota_{\cal
A}\circ\iota_H\otimes\Delta(T)=$$
$$\iota_H\otimes\delta'_K\otimes\iota_{\cal
A}(\iota_H\otimes\Delta(T)).$$
Since the norm on the range space
coincides with the norm inherited from $H\otimes{\cal A}\otimes{\cal A}$,
and since $\iota_H\otimes\Delta: H\otimes{\cal A}\to H\otimes{\cal
A}\otimes{\cal
A}$ is bounded, $\text{Ind}_0(u)$ is bounded as well.

c)
We show  relation $(3.1)$.
For 
$\psi,\psi'\in H$, $a,a'\in{\cal A}$,
$T=\psi\otimes a,T'=\psi'\otimes a'$,
$$(\iota_H\otimes\Delta(T),\iota_H\otimes\Delta(T'))_{{\cal
A}}=(\psi,\psi')h\otimes\iota(\Delta(a^*a'))=(\psi,\psi')h(a^*a')=
(T,T')I,$$
so the relation holds a fortiori on 
$X_{\text{Ind}_0(u)}$, and $\iota_H\otimes\Delta$ extends
on the completion $H_{\text{Ind}(u)}$ to a map satisfying the same
relation. By a), $\text{Ind}(u)(X_{\text{Ind}_\infty(u)}){\cal A}_\infty=
X_{\text{Ind}_\infty(u)}\odot{\cal A}_\infty$, and this subspace is
norm dense in the right Hilbert ${\cal
A}$--module $H_{\text{Ind}(u)}\otimes{\cal A}$.
\medskip

We shall call  $\text{Ind}(u)$  the {\it  representation induced}
from $u$.
We conclude this section with a result on  Frobenius reciprocity
for induced representations.
\medskip

\noindent{\bf 5.5 Theorem} {\sl 
Let $K$ be a compact quantum subgroup of the compact quantum group $G$.
Let $u$ and $v$ be finite dimensional unitary representations of $K$
and $G$ respectively, with $v$ faithful. Then the spaces
$(v,\text{Ind}_\infty(u))$
and $(v\upharpoonright_K, u)$ are linearly isomorphic.}\medskip

\noindent{\bf Proof}
Let $H_u$ and $H_v$ denote the spaces of $u$ and $v$ respectively.
For $T\in(v,\text{Ind}_\infty(u))$ and $\psi\in H_v$, $\iota\otimes e
(T(\psi))$, with $e$ the counit of ${\cal A}_\infty$ \cite{WLesHouches},
is an element of $H_u$. So we get  a  linear
map $S$ from $H_v$ to $H_u$. We show that this map is an
intertwiner between
the
desired representations. 
$$ u\circ S(\psi)=u(\iota\otimes e(T(\psi)))=
\iota_{H_u}\otimes\iota_{{\cal A}'}\otimes e(u\otimes \iota_{{\cal
A}}(T\psi))=$$
$$\iota_{H_u}\otimes\iota_{{\cal A}'}\otimes
e(\iota_{H_u}\otimes\delta'_K(T\psi))=
\iota_{H_u}\otimes\iota_{{\cal A}'}\otimes
e(\iota_{H_u}\otimes\pi\otimes\iota_{{\cal 
A}}\circ\iota_{H_u}\otimes\Delta(T\psi))=$$
$$\iota_{H_u}\otimes\pi\circ\iota_{H_u}\otimes\iota_{{\cal
A}}\otimes e\circ\iota_{H_u}\otimes\Delta(T\psi)=
\iota_{H_u}\otimes\pi
(T\psi).$$
On the other hand: 
$$S\otimes\iota_{{\cal A}'} v\upharpoonright_K(\psi)=
\iota_{H_u}\otimes e\otimes\iota_{{\cal A}'}\circ T\otimes\iota_{{\cal
A}'}\circ\iota_{H_v}\otimes\pi(v(\psi))=$$
$$\iota_{H_u}\otimes e\otimes\pi(T\otimes\iota_{{\cal
A}}(v(\psi)))=
\iota_{H_u}\otimes e\otimes\pi\circ
\iota_{H_u}\otimes\Delta
(T(\psi))=$$
$$\iota_{H_u}\otimes\pi\circ
\iota_{H_u}\otimes(e\otimes\iota_{\cal A}\circ \Delta)
T(\psi)=
\iota_{H_u}\otimes\pi(T(\psi)).$$
Let now start from an operator $S\in(v\upharpoonright_K, u)$.
For a vector $\psi\in H_v$ we 
set 
$$T(\psi):=S\otimes\iota_{{\cal A}}(v(\psi))\in H_u\odot{\cal A}_\infty.$$
We show that $T(\psi)$ lies in the space of $\text{Ind}_\infty(u)$.
$$u\otimes\iota_{{\cal A}}(S\otimes\iota_{{\cal A}}(v(\psi)))=
(uS)\otimes\iota_{{\cal A}}(v(\psi))=$$
$$(S\otimes\iota_{{\cal A}'}\circ
v\upharpoonright_K)\otimes\iota_{{\cal A}}(v(\psi))=
S\otimes\iota_{{\cal A}'}\otimes\iota_{\cal
A}\circ\iota_{H_v}\otimes\pi\otimes\iota_{\cal A}\circ
v\otimes\iota_{{\cal
A}}(v(\psi))=$$
$$S\otimes\iota_{{\cal A}'}\otimes\iota_{\cal
A}\circ\iota_{H_v}\otimes\pi\otimes\iota_{\cal A}\circ
\iota_{H_v}\otimes\Delta
(v(\psi))=
S\otimes\iota_{{\cal A}'}\otimes\iota_{\cal
A}\circ
\iota_{H_v}\otimes\delta'_K
(v(\psi))=$$
$$\iota_{H_u}\otimes\delta'_K(S\otimes\iota_{{\cal A}}(v(\psi))).$$
We check that $T\in(v,\text{Ind}_\infty(u)).$
$$\iota_{H_u}\otimes\Delta
T(\psi)=\iota_{H_u}\otimes\Delta
(S\otimes\iota_{{\cal A}}(v(\psi)))=$$
$$S\otimes\iota_{{\cal A}}\iota_{{\cal
A}}(\iota_{H_v}\otimes\Delta(v(\psi)))=
S\otimes\iota_{{\cal A}}\otimes\iota_{{\cal A}}(v\otimes\iota_{{\cal
A}}(v(\psi))=T\otimes\iota_{\cal A}(v(\psi)).$$
Finally we check that the maps $T\to S$ and $S\to T$ are inverses of one
another.
$$(\iota_{H_u}\otimes e\circ T)\otimes\iota_{{\cal A}}(v(\psi))=
\iota_{H_u}\otimes e\otimes\iota_{{\cal
A}}\circ\iota_{H_u}\otimes\Delta(T\psi)=T(\psi),$$
$$\iota_{H_u}\otimes e(S\otimes\iota_{{\cal A}}v(\psi))=
S\otimes e(v(\psi))=S(\psi),$$
since $v$ is faithful.

\end{section}

\begin{section}{Ergodicity and transitivity}

It is well known that transitivity characterizes $G$--actions arising 
from compact subgroups defined up to conjugacy, in the following
way. 
Let $X$ be a compact topological space on which a compact group $G$ acts
continuously by homeomorphisms on the right. If the $G$--action is
transitive,
the stabilizer of a point $x\in X$ is
a closed
subgroup $G_x$ of $G$, and the map $xg\in X\to G_xg\in G_x\backslash G$
is a
homeomorphism. Stabilizers of different points are conjugate
closed subgroups of $G$. Conversely, given a closed subgroup $K$ of $G$,
the quotient space $K\backslash G$ with its quotient topology becomes a
quotient
right 
$G$--space under the action $g: Kg_1\to Kg_1g$. 
The following fact is well known.
\medskip

\noindent{\bf 6.1 Proposition} {\sl Let $X$ be a compact right $G$--space
over a compact group $G$.
The following
properties are equivalent:
\begin{description}
\item{\rm a)} the automorphic action $\alpha: g\in
G\to\alpha_g\in\text{Aut}({\cal C}(X))$, with
$\alpha_g(f)(x)=f(xg)$, is ergodic,
\item{\rm b)} the $G$--action on $X$ by homeomorphisms is
transitive,
\item{\rm c)} there is no proper closed subset $F\subset X$ such that
$FG=F$.
\end{description}
If one of the above properties holds then
 there is a closed subgroup $K$ of $G$, determined up to conjugation, such
that
 (${\cal C}(X)$, $\alpha$) is isomorphic to the automorphic
action on ${\cal C}(K\backslash G)$ induced by the natural right
$G$--action on
$K\backslash G$.
}
\medskip

Let now $G$ be a compact quantum group acting on a (possibly noncommutative)
unital $C^*$--algebra ${\cal C}$.
Is ergodicity still equivalent to some sort of topological transitivity of
that
action?

In the case where 
 $G$  is a classical compact group
acting pointwise continuously in norm on a  unital
$C^*$--algebra ${\cal
C}$,
$\alpha:G\to\text{Aut}({\cal C})$, the second question has been
investigated by Longo and Peligrad in \cite{LP}, who proved the following
theorem.

Recall that a subset $M$ of ${\cal C}$ is called
$G$--invariant if $\alpha_g(M)\subset M$ for all
$g\in G$. Clearly, since each $g$ has an inverse $g^{-1}$,
$\alpha_g(M)=M$ for all $g\in G$.
\medskip

\noindent{\bf 6.2 Theorem} {\sl {\rm \cite{LP}} The following properties
are
equivalent:
\begin{description}
\item{\rm a)} ${\cal C}^\alpha={\mathbb C}I$,
\item{\rm b)} there is no proper, closed, $G$--invariant left ideal 
of ${\cal C}$,
\item {\rm c)} there is no proper, hereditary,  $G$--invariant 
$C^*$--subalgebra of ${\cal C}$.
\end{description}
}\medskip

We give a generalization of the above result in the quantum group case.
A $C^*$--subalgebra, or a closed left ideal,  ${\cal B}$ of ${\cal C}$
will be 
called $G$--invariant if $\delta({\cal B})\subset{\cal B}\otimes{\cal
A}$.
\medskip

\noindent{\bf 6.3 Theorem} {\sl Let $\delta:{\cal C}\to{\cal
C}\otimes{\cal
A}$ an action of the compact quantum group $G=({\cal A},\Delta)$
on the unital $C^*$--algebra ${\cal C}$ such that $I\otimes{\cal
A}\delta({\cal C})={\cal C}\otimes{\cal A}$.
Then the following properties are equivalent:
\begin{description}
\item{\rm a)} ${\cal C}^\delta={\mathbb C}I$,
\item{\rm b)} there is no proper, closed $G$--invariant, left ideal ${\cal
I}$ of ${\cal C}$ such that $I\otimes{\cal A}\delta({\cal I})$ is
norm dense
in ${\cal I}\otimes{\cal A}$,
\item{\rm c)} there is no proper, $G$--invariant, hereditary
$C^*$--subalgebra ${\cal H}$ of ${\cal C}$ such that
the hereditary $C^*$--subalgebra of ${\cal C}\otimes{\cal A}$
generated by $\delta({\cal H})$ is ${\cal H}\otimes{\cal A}$,
\item{\rm d)} there is no proper, open projection $p\in{\cal C}''$ such
that
$\delta''(p)=p\otimes I$.
\end{description}}\medskip

\noindent{\bf Proof} 
Closed left ideals of ${\cal C}$ are in bijective correspondence
with hereditary $C^*$--subalgebras via the map assigning to the left ideal
${\cal J}$ the algebra ${\cal H}={\cal J}\cap{\cal J}^*$ with
inverse the map assigning to the algebra ${\cal H}$ the ideal
${\cal J}=\{j\in{\cal C}: j^*j\in{\cal H}\}$ (see \cite{Pedersen}).
If ${\cal J}$ and ${\cal H}$ are in correspondence, 
for any $C^*$--algebra ${\cal A}$, ${\cal J}\otimes{\cal A}$
and ${\cal H}\otimes{\cal A}$ are a closed left ideal and a hereditary
$C^*$--subalgebra of ${\cal C}\otimes{\cal A}$ in correspondence as well.
If ${\cal J}$ is $G$--invariant and $h\in{\cal H}$ then
$\delta(h)\in{\cal J}\otimes{\cal A}\cap{\cal J}^*\otimes{\cal A}={\cal
H}\otimes{\cal A}$, so ${\cal H}$ is $G$--invariant.
Conversely, if ${\cal H}$ is $G$--invariant and $j\in{\cal J}$
then $\delta(j^*j)\in{\cal H}\otimes{\cal A}$. Therefore
$\delta(j)\in{\cal J}\otimes{\cal A}$. Thus $G$--invariant
closed left ideals are in bijective correspondence with $G$--invariant
hereditary $C^*$--subalgebras.
We claim that if ${\cal I}$ is a closed left ideal
of ${\cal C}$ then the closure ${\cal J}$ of $I\otimes{\cal A}\delta({\cal
I})$ is
a left  ideal of ${\cal C}\otimes{\cal A}$. We only have to show that 
${\cal C}\otimes I{\cal J}\subset{\cal J}$. Now since $I\otimes{\cal
A}\delta({\cal C})$ is dense, ${\cal C}\otimes I$ is contained in the 
closure of ${I}\otimes{\cal A}\delta({\cal C})$, so ${\cal C}\otimes
I{\cal J}$, which coincides with the closure of $I\otimes{\cal
A}{\cal C}\otimes I\delta({\cal I})$,
is contained in the closure of $I\otimes{\cal A}\delta({\cal C}{\cal I})$,
which is in turn contained in ${\cal J}$. 
By \cite{Pedersen}, hereditary $C^*$--subalgebras, or closed left
ideals,
of a unital $C^*$--algebra ${\cal B}$ are in bijective correspondence
with open projections of ${\cal B}''$. If $p$ is the projection
of ${\cal C}''$ corresponding to a closed, left, $G$--invariant ideal
${\cal I}$
of ${\cal C}$ then $\delta''(p)$ must coincide with the open
projection of $({\cal C}\otimes{\cal A})''$ corresponding to ${\cal J}$,
the closure of the 
ideal $I\otimes{\cal A}\delta({\cal I})$;
$G$--invariance of ${\cal I}$ shows that this ideal is contained
in ${\cal I}\otimes{\cal A}$, so
$\delta''(p)
\leq
p\otimes
I$.
If in addition $I\otimes{\cal A}\delta({\cal I})$ is norm dense
in ${\cal I}\otimes{\cal A}$ then $p\otimes I\delta''(p)=p\otimes I$,
so $\delta''(p)=p\otimes I$. Conversely, any such projection corresponds
to a closed $G$--invariant left ideal ${\cal I}$ such that $I\otimes{\cal
A}\delta({\cal I})={\cal I}\otimes{\cal A}$.
Again, for a hereditary $G$--invariant $C^*$--subalgebra ${\cal H}$
of ${\cal C}$, the hereditary $C^*$--subalgebra generated by $\delta({\cal
H})$ corresponds to the projection $\delta(p)$, which is in turn dominated
by $p\otimes I$ by $G$--invariance. Requiring that it
coincides 
with ${\cal H}\otimes{\cal A}$ fixes $\delta''(p)=p\otimes I$.
We have thus proven the equivalence of b), c) and d).
The implication b)$\to$ a) is easy: if
 there were a positive non scalar element $a$ in ${\cal C}^\delta$
then its spectrum would contain at least two points $x_1$ and $x_2$.
Let $f$ be a continuous function on the spectrum such that $f(x_1)=1$,
$f(x_2)=0$. The closed left ideal ${\cal J}$ generated by $f(a)$ is then
proper and
$G$--invariant. The norm closure of $I\otimes{\cal A}\delta({\cal
J})$ is genereated by $I\otimes{\cal A}\delta({\cal C})f(a)\otimes I$,
or by
${\cal C}f(a)\otimes{\cal A}$ since $I\otimes{\cal A}\delta({\cal C})$ is
dense
in ${\cal C}\otimes{\cal A}$. We are left to show that a)$\to$d).
Let $p$ be an open projection of ${\cal C}''$ such that
$\delta''(p)=p\otimes I$. Consider the conditional expectation 
$E:{\cal C}\to{\cal C}^\delta$ over the fixed points obtained averaging
over $G$: $E(c):=\iota\otimes h\circ\delta(c)$. 
By universality of ${\cal C}''$,
$E$ extends to a normal positive map $E'':{\cal C}''\to{\cal C}''$
such that $E''(c'')=c''$ whenever $\delta''(c'')=c''\otimes I$.
In particular, $E''(p)=p$. By 3.11.9 in \cite{Pedersen}
$p$ can
be obtained as a strong
limit of a bounded monotone increasing net $x_\alpha$ from ${\cal
C}^+$, therefore
$E''(p)=p$ is the strong limit of the monotone increasing net
$E(x_\alpha)$.
Since 
${\cal C}^\delta={\mathbb C}I$, $p$ must be a multiple of the identity,
i.e. either $p=0$, $p=I$, and the proof is complete.
\medskip

Under stronger assumptions 
on the quantum group $G$, conditions b) through d) in the previous theorem
take a more relaxed form.\medskip

\noindent{\bf 6.4  Theorem} {\sl If  the Haar measure $h$ and
the action $\delta$ are faithful maps, then, under the same assumption
as in the previous theorem, the following conditions are equivalent:
\begin{description}
\item{\rm a)} ${\cal C}^\delta={\mathbb C}$,
\item{\rm b')} there is no proper closed $G$--invariant left ideal of
${\cal
C}$,
\item{\rm c')} there is no proper hereditary $G$--invariant
$C^*$--subalgebra of ${\cal C}$,
\item{\rm d')} there is no proper open projection $p\in{\cal C}''$ such
that $\delta(p)\leq p\otimes I$.
\end{description}}\medskip

\noindent{\bf Proof}
Indeed, projections as in d') are in one to one correspondence with ideals
as in b') and algebras as in c'). If $p\in{\cal C}''$ is a
nonzero projection satisfying the condition 
stated in d') then $E''(p)\leq p$ and $p$ is the strong limit of 
an increasing net $x_\alpha$ from ${\cal C}^+$ with nonzero elements,
so  $E''(p)$ turns out to be a strong limit of the increasing net
$E(x_\alpha)$, with nonzero entries. 
If the action is ergodic then $E''(p)$ is a nonzero scalar,
and therefore  $p=I$.
\medskip

\end{section}

\begin{section}{Ergodicity and  quotient spaces} 

When
is an ergodic action $\delta:{\cal C}\to{\cal C}\otimes{\cal A}$ of a 
compact quantum group $G=({\cal A}, \Delta)$ on a 
$C^*$--algebra, isomorphic
to an action
 on some quantum quotient space $K\backslash G=({\cal A}^{\delta'_K},
\eta')$?

Ergodicity is clearly not enough. Just think  
of the case where $G$ is a group. The question has a negative answer
if ${\cal C}$ is not commutative, and positive
 if  ${\cal C}$ is commutative, thanks to Prop. 6.1.

But, of course, commutativity of ${\cal C}$ can not be assumed in the quantum
case. We look for some other 
property, necessary also in the case where $G$ is a quantum
group.
The following is a key observation.
\medskip

\noindent{\bf 7.1 Proposition} {\sl Let $\alpha:G\to\text{Aut}({\cal C})$
be a strongly continuous, ergodic, action of a compact group on a unital
$C^*$--algebra ${\cal C}$. Then the following conditions are equivalent:
\begin{description}
\item{\rm a)} ${\cal C}$ has a character,
\item{\rm b)} ${\cal C}$ is commutative.
\end{description}
If one of the above condition is satisfied, the given dynamical
system is isomorphic to the system arising from a quotient
right
$G$--space.
}\medskip

\noindent{\bf Proof}
We need to show that a)$\to$b). 
Let $\chi$ be a character of
${\cal C}$. 
The  stabilizer of $\chi$, $G_\chi:=\{ g\in G:
\chi\circ\alpha_g=\chi\}$, is  a closed subgroup of $G$,
and the map
$$\rho: c\in{\cal C}\to (g\in G\to\chi\circ\alpha_{g}(c))\in {\cal C}(G)$$
is a $^*$--homomorphism  
with range included in the commutative $C^*$--algebra
${\cal C}(G_\chi\backslash
G)$.
This map intertwines the corresponding automorphic
$G$--actions. 
The two sided closed ideal of ${\cal C}$:
${\cal J}=\{c\in{\cal C}: \chi(\alpha_g(c))=0, g\in G\}$ is
obviously $G$--invariant and does not contain the identity, so
by ergodicity, ${\cal J}=0$.
It follows that the $^*$--homomorphism ${\cal C}\to {\cal C}(G)$ assigning
to an
element
$c\in{\cal C}$
the continuous function $g\in G\to\chi(\alpha_g(c))$,
is faithful. Thus ${\cal C}$ is commutative, and $G$ acts transitively on
the spectrum
of ${\cal C}$. In particular, the closed orbit
$\{\chi\circ\alpha_g, g\in G\}$ must coincide with whole spectrum of
${\cal C}$. 
A straightforward application of the Stone--Weierstrass theorem shows that
$\rho$ must be surjective.
\medskip

How many  quantum coset spaces with a $^*$--character do there exist?
Think of the following construction. A quantum coset space
is the fixed point algebra of the Hopf $C^*$--algebra ${\cal A}$
of a quantum group $G=({\cal A},\Delta)$
under the action of a subgroup  (see Sect. 4). Therefore it 
suffices to look for a $^*$--character on ${\cal A}$.

On the other hand Woronowicz shows in \cite{Wcmp} that
every compact matrix pseudogroup has a densely defined
$^*$--character $e$: the counit. 
This is a $^*$--homomorphism $e:{\cal A}\to{\mathbb C}$, defined
only on the smooth part ${\cal A}_\infty$ of ${\cal A}$,
such that for $a\in{\cal A}_\infty$,
$$\iota\otimes e\circ\Delta(a)=a.$$
But, he also shows in \cite{Wtk} that one can obtain
 compact quantum groups
 from any suitable category of
finite dimensional Hilbert spaces, via Tannaka--Krein duality. 
In fact these groups obtained from categories are  completion of their
smooth part with respect to the maximal $C^*$--norm. Therefore
for these groups the counit must be
an everywhere defined $^*$--character.
For example, the group $S_\mu U(d)$ has such a character.

In conclusion, for sufficiently many quantum groups $G$,
if a compact quantum subgroup $K$ of $G$ is given, the restriction $e_K$
 of the counit $e$ to the right coset space ${\cal A}^{\delta'_K}$
is a continuous $^*$--character of that $C^*$--subalgebra of ${\cal A}$.

We show a property possessed by the counit.
\medskip

\noindent{\bf 7.2 Lemma} {\sl If the action $\delta:{\cal C}\to{\cal
C}\otimes{\cal A}$ of the compact quantum group $G$ on ${\cal A}$
is faithful, and if $e$ is an everywhere defined counit
of $G$ then $\iota\otimes e\circ\delta(c)=c$ for all $c\in{\cal C}$.}
\medskip

\noindent{\bf Proof}
It suffices to show that 
$\delta(\iota\otimes e\circ\delta(c))=\delta(c)$.
Indeed, the l.h.s. equals
$$\iota\otimes\iota\otimes e\circ\delta\otimes\iota\circ\delta(c)=
\iota\otimes\iota\otimes
e\circ\iota\otimes\Delta\circ\delta(c)=$$
$$\iota\otimes(\iota\otimes e\circ\Delta)(\delta(c))=\delta(c).$$
\medskip

Consider an action $\delta$ of a compact quantum group $G=({\cal
A},\Delta)$
on a unital $C^*$--algebra ${\cal C}$. If ${\cal C}$ has a character
$\chi$, Prop. 7.1 suggests how to construct an intertwiner from the system
$({\cal C}, \delta,\chi)$ to some system of the form $K\backslash G=({\cal
A}^{\delta'_K}, \eta', e_K)$.
\medskip

\noindent{\bf 7.3 Theorem} {\sl 
Let $\delta:{\cal C}\to{\cal C}\otimes{\cal A}$ be an action
of $G$ on a unital $C^*$--algebra ${\cal C}$.
Let  $\chi:{\cal C}\to{\mathbb C}$ be a character of ${\cal C}$.
Then there is a compact quantum subgroup $G_\chi=({\cal A}',\Delta')$ of
$G$ 
such that the map
$$T_\chi:=\chi\otimes\iota_{\cal A}\circ\delta:{\cal C}\to{\cal
A}^{\delta'_{G_\chi}}$$
is a $^*$--homomorphism intertwining the corresponding $G$--actions:
$$\eta'\circ T_\chi=T_\chi\otimes\iota_{\cal A}\circ\delta.$$
Furthermore, if $\pi:{\cal A}\to{\cal A}'$ is the quotient map,
$$\pi(T_\chi(c))=\chi(c) I_{{\cal A}'},\quad c\in {\cal C}.$$
Also, if $\delta$ is faithful and if $G$ has an everywhere defined counit
$e$ then
$$e(T_\chi(c))=\chi(c),\quad c\in{\cal C}.$$
}
\medskip

\noindent{\bf Proof}
Consider the linear subspace $M_\chi$ of ${\cal A}$ 
generated by
$$\{\chi\otimes\iota\circ\delta(c)-\chi(c)I, c\in{\cal C}\}.$$
Notice that
$$\Delta(\chi\otimes\iota\circ\delta(c)-\chi(c)I)=
(\chi\otimes\iota\circ\delta-\chi
I)\otimes\iota(\delta(c))+I\otimes(\chi\otimes\iota\circ\delta(c)-\chi(c)I),$$
so 
$$\Delta(M_\chi)\subset M_\chi\otimes{\cal A}+I\otimes M_\chi.$$
Let $J_\chi$ be the closed ideal generated by $M_\chi$ in ${\cal A}$.
This is a $^*$--ideal since $M_\chi$ is $^*$--invariant.
The above relation shows that
$\Delta(J_\chi)$ is contained in 
$J_\chi\otimes
{\cal A}+{\cal A}\otimes J_\chi.$
Set ${\cal A}'={\cal
A}/J_\chi$, and define
$\Delta'([a])=\pi\otimes \pi(\Delta(a))$,
with $\pi:{\cal A}\to{\cal A}'$ the quotient map. 
This map is well defined and defines a nondegenerate coassociative coproduct 
on ${\cal A'}$.
It is now obvious that $G_\chi:=({\cal A}', \Delta')$ is a compact
subgroup of $G$, and, by definition,
$$\chi\otimes\pi\circ\delta(c)=\chi(c)I_{{\cal A}'},\quad c\in{\cal C}.$$
We show that the range of $T_\chi$ is included
in ${\cal A}^{\delta'_{G_\chi}}$. For $c\in{\cal C}$,
$$\delta'_{G_\chi}(T_\chi(c))=
\pi\otimes\iota\circ\Delta\circ\chi\otimes\iota\circ\delta(c)=$$
$$\chi\otimes\pi\otimes\iota\circ\iota\otimes\Delta\circ\delta(c)=
\chi\otimes\pi\otimes\iota\circ\delta\otimes\iota\circ\delta(c)=$$
$$[(\chi\otimes\pi\circ\delta)\otimes\iota]\circ\delta(c)=\chi(\ 
.\ )I\otimes\iota\circ\delta(c)=I\otimes T_\chi(c).$$
We finally show that $T_\chi$ intertwines the corresponding
$G$--actions. For $c\in{\cal C}$,
$$\eta'\circ T_\chi(c)=\Delta\circ\chi\otimes\iota\circ\delta(c)=$$
$$\chi\otimes\iota\otimes\iota\circ\iota\otimes\Delta\circ\delta(c)=
\chi\otimes\iota\otimes\iota\circ\delta\otimes\iota\circ\delta(c)=$$
$$T_\chi\otimes\iota\circ\delta(c).$$
Finally, by the previous lemma,
$$e(T_\chi(c))=\chi\otimes e\circ\delta(c)=\chi(\iota\otimes
e\circ\delta(c))=\chi(c).$$
\medskip

The subgroup $G_\chi$  constructed along the proof of the previous theorem
will be called the {\it  subgroup stabilizing   $\chi$}.

We can  summarize the results of this and the previous section in
the following
theorem.
\medskip

\noindent{\bf 7.4 Theorem} {\sl 
Let $\delta:{\cal C}\to{\cal C}\otimes{\cal A}$ be an action 
of the compact quantum group $G=({\cal A},\Delta)$ on a unital 
$C^*$--algebra ${\cal C}$ endowed with a $^*$--character $\chi$. Assume
that 
\begin{description}
\item{\rm a)} $I\otimes{\cal A}\delta({\cal C})$ is norm dense in ${\cal
C}\otimes{\cal A}$,
\item{\rm b)} the action $\delta$ is ergodic: ${\cal C}^\delta={\mathbb
C}$.
\end{description}
Assume furthermore that
 the action $\delta$ and the Haar measure $h$ of $G$
are faithful maps.
It follows that $T_\chi$ is faithful. 
}\medskip

\noindent{\bf Proof}
We show that the kernel ${\cal J}$ of $T_\chi$ is $G$--invariant.
For $j\in{\cal J}$, by Theorem 7.3, $T_\chi\otimes\iota_{\cal
A}(\delta(j))=\eta'(T_\chi(j))=0$, so $\delta(j)\in\text{ker
}(T_\chi\otimes\iota_{\cal A})={\cal J}\otimes{\cal A}$. Since ${\cal
J}\neq{\cal C}$, we must have ${\cal J}=0$ thanks to Theorem 6.4.
\medskip

\noindent{\bf 7.5 Corollary} {\sl If $\delta:{\cal C}\to{\cal
C}\otimes{\cal
A}$ is a faithful ergodic action of a compact quantum group $G=({\cal
A},\Delta)$, with faithful Haar measure,
on a commutative unital $C^*$--algebra ${\cal C}$ satisfying the
nondegeneracy property a) above, then $({\cal C},\delta)$ can be embedded
faithfully into a quotient space of $G$.}\medskip

\noindent{\it Remark 6} If we drop the assumption that $h$ and $\delta$
are
faithful maps, but if we know a priori that $T_\chi$ is surjective, then
$T_\chi$ must be faithful, and therefore  a $^*$--isomorphism. 
The reason is explained in the following lemma, which we include for
future reference.
\medskip

\noindent{\bf 7.6 Lemma} {\sl Let $\delta:{\cal C}\to{\cal C}\otimes{\cal
A}$
and $\eta:{\cal D}\to{\cal D}\otimes{\cal A}$ be unital actions of the
compact quantum group $G=({\cal A},\Delta)$, and let $T:{\cal C}\to{\cal
D}$ be a surjective  $^*$--homomorphism such that
$$\eta\circ T=T\otimes\iota_{\cal A}\circ \delta.$$
If $I\otimes{\cal A}\delta({\cal C})$ and $I\otimes{\cal
A}\eta({\cal D})$ are norm dense in ${\cal C}\otimes{\cal A}$ and  ${\cal
D}\otimes{\cal
A}$ respectively, then the kernel ${\cal J}$ of $T$ is a closed, two
sided,
$G$--invariant ideal of ${\cal C}$ such that 
$I\otimes{\cal A}\delta({\cal
J})$ is norm dense in ${\cal J}\otimes{\cal A}$.
In particular, if $\delta$ is ergodic, $T$ must be faithful.
}\medskip

\noindent{\bf Proof} 
As in the proof of the previous theorem, one can show that ${\cal J}$ is 
$G$--invariant. By Theorem 6.3,
we are  left to show that $I\otimes{\cal A}\delta({\cal J})$
is dense in ${\cal J}\otimes{\cal A}$. Extend $T$ to a normal
$^*$--homomorphism 
$T'':{\cal C}''\to{\cal D}''$, and let $p$ be the central projection
of ${\cal C}''$ such that $\text{ker}(T'')={\cal C}''(I-p)$.
$I-p$ is the open projection of ${\cal C}''$ corresponding to
${\cal J}$. We need to show, by Theorem 6.3, that
$\delta''(I-p)=(I-p)\otimes I$.
The restriction $T_p$ of $T''$ to ${\cal C}''p$ is a
normal $^*$--monomorphism with range ${\cal D}''$. Since 
$I\otimes{\cal A}\eta({\cal D})$ is norm dense in ${\cal D}\otimes{\cal
A}$, it is a fortiori weakly dense in the von Neumann tensor product
${\cal D}''\otimes{\cal A}''$. Pulling back this relation with  
$T_p^{-1}\otimes{ I}$ shows that 
$I\otimes{\cal A}''\delta''({\cal C}''p)$ 
is weakly dense in the
von Neumann tensor product ${\cal
C}''p\otimes{\cal A}''$. Now the property that $I\otimes{\cal
A}\delta({\cal C})$ is norm dense shows that
the weak closure of 
$I\otimes{\cal A}''\delta''({\cal C}''p)$ 
is a weakly closed left (and right) ideal of ${\cal C}''\otimes{\cal A}''$
defined by the projection $\delta''(p)$. The density statement
shows that $\delta''(p)=p\otimes I$, and therefore
$\delta''(I-p)=(I-p)\otimes I$.
\medskip  

\noindent{\it Remark 7}
We conclude the paper noting that one can not expect in general 
an isomorphism of an ergodic $G$--space $({\cal C},\delta)$ with 
a quotient $G$--space $K\backslash G$ by a stabilizer subgroup.
In fact, Wang shows in \cite{Wang} an example of a compact quantum group
acting ergodically on a commutative $C^*$--algebra for which 
the quotient space by a point stabilizer subgroup is not commutative.
\bigskip

\noindent{\bf Acknowledgments}
I would like to warmly thank Shuzhou Wang for his interest in this paper, and for kindly communicating to me 
  a gap in Lemma 2.3 of the published version, which resulted in the current revised version,   which also benefits from a   comparison between  the  notion of a closed coideal in the sense of this paper,   the notion of a Woronowicz 
$C^*$--ideal  \cite{Wang_free} and issues on density of the smooth algebraic subcoideal.

\end{section}

\end{document}